\documentclass[reqno,11pt]{amsart}
\usepackage{amscd,amssymb}

\theoremstyle{plain}
\newtheorem{Thm}[subsection]{Theorem}
\newtheorem{Cor}[subsection]{Corollary}
\newtheorem{Lem}[subsection]{Lemma}
\newtheorem{Prop}[subsection]{Proposition}
\newtheorem{Conj}[subsection]{Conjecture}

\theoremstyle{definition}
\newtheorem{Def}[subsection]{Definition}

\theoremstyle{remark}

\newtheorem{Rem}[subsection]{Remark}

\newtheorem{conj}{Conjecture}

\numberwithin{equation}{section}
\renewcommand{\rm}{\normalshape}

\newif\ifShowLabels
\ShowLabelstrue
\newdimen\theight
\def\TeXref#1{%
    \leavevmode\vadjust{\setbox0=\hbox{{\tt
        \quad\quad  {\small \rm #1}}}%
    \theight=\ht0
    \advance\theight by \lineskip
    \kern -\theight \vbox to
    \theight{\rightline{\rlap{\box0}}%
    \vss}%
    }}%

\ShowLabelsfalse

\renewcommand{\sec}[2]{\section{#2}\label{S:#1}%
    \ifShowLabels \TeXref{{S:#1}} \fi}
\newcommand{\ssec}[2]{\subsection{#2}\label{SS:#1}%
    \ifShowLabels \TeXref{{SS:#1}} \fi}

\newcommand{\refs}[1]{Section ~\ref{S:#1}}
\newcommand{\refss}[1]{Section ~\ref{SS:#1}}

\newcommand{\reft}[1]{Theorem ~\ref{T:#1}}
\newcommand{\refl}[1]{Lemma ~\ref{L:#1}}
\newcommand{\refp}[1]{Proposition ~\ref{P:#1}}

\newcommand{\refd}[1]{Definition ~\ref{D:#1}}

\newcommand{\refe}[1]{\eqref{E:#1}}
\newcommand{\refco}[1]{Conjecture ~\ref{Co:#1}}

\newenvironment{thm}[1]%
    { \begin{Thm} \label{T:#1}  \ifShowLabels \TeXref{T:#1} \fi }%
    { \end{Thm} }

\renewcommand{\th}[1]{\begin{thm}{#1} \sl }
\renewcommand{\eth}{\end{thm} }

\newenvironment{lemma}[1]%
    { \begin{Lem} \label{L:#1}  \ifShowLabels \TeXref{L:#1} \fi }%
    { \end{Lem} }
\newcommand{\lem}[1]{\begin{lemma}{#1} \sl}
\newcommand{\elem}{\end{lemma}}

\newenvironment{propos}[1]%
    { \begin{Prop} \label{P:#1}  \ifShowLabels \TeXref{P:#1} \fi }%
    { \end{Prop} }
\newcommand{\prop}[1]{\begin{propos}{#1}\sl }
\newcommand{\eprop}{\end{propos}}

\newenvironment{corol}[1]%
    { \begin{Cor} \label{C:#1}  \ifShowLabels \TeXref{C:#1} \fi }%
    { \end{Cor} }
\newcommand{\cor}[1]{\begin{corol}{#1} \sl }
\newcommand{\ecor}{\end{corol}}

\newenvironment{defeni}[1]%
    { \begin{Def} \label{D:#1}  \ifShowLabels \TeXref{D:#1} \fi }%
    { \end{Def} }
\newcommand{\defe}[1]{\begin{defeni}{#1} \sl }
\newcommand{\edefe}{\end{defeni}}

\newenvironment{remark}[1]%
    { \begin{Rem} \label{R:#1}  \ifShowLabels \TeXref{R:#1} \fi }%
    { \end{Rem} }
\newcommand{\rem}[1]{\begin{remark}{#1}}
\newcommand{\erem}{\end{remark}}

\newenvironment{conjec}[1]%
    { \begin{Conj} \label{Co:#1}  \ifShowLabels \TeXref{Co:#1} \fi }%
    { \end{Conj} }
\renewcommand{\conj}[1]{\begin{conjec}{#1} \sl }
\newcommand{\econj}{\end{conjec}}

\newcommand{\eq}[1]%
    { \ifShowLabels \TeXref{E:#1} \fi
       \begin{equation} \label{E:#1} }
\newcommand{\eeq}{ \end{equation} }

\newcommand{\prf}{ \begin{proof} }
\newcommand{\epr}{ \end{proof} }

\newcommand\alp{\alpha}     

     \newcommand\Gam{\Gamma}
\newcommand\del{\delta}     
\newcommand\eps{\varepsilon}

\newcommand\lam{\lambda}        \newcommand\Lam{\Lambda}
\newcommand\sig{\sigma}

\newcommand\calF{{\mathcal{F}}}

\newcommand\calH{{\mathcal{H}}}

\newcommand\calK{{\mathcal{K}}}
\newcommand\calL{{\mathcal{L}}}

\newcommand\calO{{\mathcal{O}}}

\newcommand\calS{{\mathcal{S}}}

\newcommand\calU{{\mathcal{U}}}

\newcommand\calW{{\mathcal{W}}}

        \newcommand\bfC{{\mathbf C}}
        \newcommand\bfD{{\mathbf D}}

        \newcommand\bfS{{\mathbf S}}

\newcommand\PP{\mathbb{P}}
\renewcommand\AA{\mathbb{A}}

\newcommand\GG{\mathbb{G}}

\newcommand\ZZ{\mathbb{Z}}

\newcommand\CC{\mathbb{C}}

\newcommand\NN{\mathbb{N}}

 \newcommand\gra{{\mathfrak{a}}}

 \newcommand\grg{{\mathfrak{g}}}

 \newcommand\grn{{\mathfrak{n}}}

 \newcommand\grt{{\mathfrak{t}}}

\newcommand\sdp{\times \hskip -0.3em {\raise 0.3ex
\hbox{$\scriptscriptstyle |$}}} 

\newcommand\GL{\operatorname{GL}}
\newcommand\Gr{\operatorname{Gr}}

\newcommand\Hom{\operatorname {Hom}}

\newcommand\Int{\operatorname{Int}}

\newcommand\Perv{\operatorname{Perv}}

\newcommand\SL{{\rm SL}}

\newcommand\Spec{\operatorname{Spec}}

\newcommand\Sym{\operatorname{Sym}}

\newcommand\oS{{\overline{S}}}

\newcommand\oeta{{\overline{\eta}}}

\newcommand\olam{{\overline{\lambda}}}

\newcommand\omu{{\overline{\mu}}}
\newcommand\onu{{\overline{\nu}}}

\newcommand\hatG{{\widehat{G}}}

\newcommand\tilV{{\widetilde{V}}}

\newcommand\x{\times}
\newcommand\ten{\otimes}

\newcommand{\ra}{\rangle}
\newcommand{\la}{\langle}

\newcommand\gd{{\bfG}^{\vee}}

\renewcommand\Spec{\operatorname{Spec}}

\newcommand\nc{\newcommand}
\nc{\opn}{\operatorname}
\renewcommand{\gd}{\grg^{\vee}}

\let\ggg=\gg
\renewcommand{\gg}{\operatorname{Gr}_G}
\newcommand{\go}{G(\calO)}

\newcommand{\ogg}{{\overline \gg}}
\newcommand{\ogl}{\ogg^{\lam}}
\newcommand{\pgg}{\text{Perv}_{\go}(\gg)}
\newcommand{\IC}{{\operatorname{IC}}}

\newcommand\gr{\operatorname{gr}}

\newcommand{\BZ}{{\mathbb Z}}

\newcommand{\BN}{{\mathbb N}}

\newcommand{\BA}{{\mathbb A}}
\newcommand{\BBA}{\overline{\mathbb A}{}}

\newcommand{\CF}{{\mathcal F}}

\newcommand{\Maps}{\overset{\circ}{Ufl}{}}

\nc\aff{\operatorname{aff}}
\nc\oGr{\overline{\Gr}}
\nc\Bun{\operatorname{Bun}}
\nc\hgrg{\widehat{\grg}}
\renewcommand\Int{\operatorname{Int}}
\nc\bInt{\overline{\Int}}
\nc\hatLam{\widehat{\Lam}}
\nc\bmu{\overline{\mu}}
\nc\bnu{\overline{\nu}}
\nc\blambda{\overline{\lam}}
\renewcommand\SL{\operatorname{SL}}
\nc\ocalW{\overline{\calW}}
\nc\pos{\operatorname{pos}}
\nc\IH{\operatorname{IH}}
\nc\fsl{\mathfrak{sl}}
\nc\fgl{\mathfrak{gl}}
\nc\Rep{\operatorname{Rep}}
\nc\Gal{\operatorname{Gal}}
\nc{\tilGr}{\widetilde{\Gr}}
\renewcommand\Maps{\operatorname{Maps}}
\nc\Pic{\operatorname{Pic}}
\nc\hgl{\widehat{\fgl}}
\nc\hsl{\widehat{\fsl}}
\dedicatory{To the memory of Izrail Moiseevich Gelfand}
\begin{document}
\title{Pursuing the double affine Grassmannian I: transversal
slices via instantons on $A_k$-singularities}
\author{Alexander Braverman and Michael Finkelberg}

\begin{abstract}This paper is the first in a series that describe a conjectural analog
of the geometric Satake isomorphism for an affine Kac-Moody group (in this paper for simplicity we
consider only untwisted and simply connected case). The usual geometric Satake isomorphism
for a reductive group $G$ identifies the tensor category $\Rep(G^{\vee})$ of finite-dimensional representations
of the Langlands dual group $G^{\vee}$ with the tensor
category $\Perv_{\go}(\gg)$ of $\go$-equivariant perverse sheaves on the affine
Grassmannian $\gg=G(\calK)/G(\calO)$ of $G$ (here $\calK=\CC((t))$ and $\calO=\CC[[t]]$). As a byproduct one gets a description of the irreducible $\go$-equivariant
intersection cohomology sheaves of the closures of $\go$-orbits in $\gg$
in terms of $q$-analogs of the weight multiplicity for finite dimensional representations
of $G^{\vee}$.

The purpose of this paper is to try to generalize the above results to the case when $G$ is replaced by the corresponding
affine Kac-Moody group $G_{\aff}$ (we shall refer to the (not yet constructed) affine Grassmannian of $G_{\aff}$ as the
{\it double affine Grassmannian}). More precisely, in this paper we construct certain varieties that should be thought of
as transversal slices to various $G_{\aff}(\calO)$-orbits inside the closure of another
$G_{\aff}(\calO)$-orbit in $\Gr_{G_{\aff}}$. We present a conjecture that computes the IC sheaf of these varieties in terms
of the corresponding $q$-analog of the weight multiplicity for the {\it Langlands dual affine group} $G_{\aff}^{\vee}$
and we check this conjecture in a number of cases.

Some further constructions (such as convolution of the corresponding perverse sheaves,
analog of the Beilinson-Drinfeld Grassmannian
etc.) will be addressed in another publication.
\end{abstract}
\maketitle
\sec{int}{Introduction}

\ssec{}{Langlands duality and the Satake isomorphism}
Let $F$ be a global field and let $\AA_F$ denote its ring of adeles.
Let $G$ be split reductive group over $F$. The classical {\em Langlands duality} predicts that irreducible automorphic
representations of $G(\AA_F)$ are closely related to the homomorphisms from the absolute Galois group $\Gal_F$ of $F$ to
the {\em Langlands dual group} $G^{\vee}$. Similarly, if $G$ is a split reductive group over a local-nonarchimedian
field $\calK$, Langlands duality predicts a relation between irredicible representations of $G(\calK)$ and homomorphisms
from $\Gal_{\calK}$ to $G^{\vee}$.

The starting point for Langlands duality is the following {\em Satake isomoprhism}. Let $\calO\subset\calK$ denote
the ring of integers of $\calK$. Then the group $G(\calK)$ is a locally compact topological group and
$G(\calO)$ is its maximal compact subgroup. One may study the {\em spherical Hecke algebra} $\calH$ of
$G(\calO)$-biinvariant compactly supported $\CC$-valued measures on $G(\calK)$. The Satake isomorphism is a canonical
isomorphism between $\calH$ and the complexified
Grothendieck ring $K_0(\Rep(G^{\vee}))$
of finite-dimensional representations of $G^{\vee}$.
\ssec{}{The geometric version} The classical Langlands duality has its geometric counterpart, usually referred to as
the {\em geometric Langlands duality}. It is based on the following geometric version of the Satake isomorphism.
Let now $\calK=\CC((s))$ and let $\calO=\CC[[s]]$;
here $s$ is a formal variable. Let $\Gr_G=G(\calK)/G(\calO)$.
Then the geometric analog of the algebra $\calH$ considered above is the category $\Perv_{G(\calO)}(\Gr_G)$ of
$G(\calO)$-equivariant perverse sheaves on $\Gr_G$ (cf. \cite{Gi}, \cite{Lu-qan} or \cite{MV}
for the precise definitions). According to {\em loc. cit.} the category $\Perv_{G(\calO)}(\Gr_G)$
possesses canonical tensor structure and the geometric Satake isomoprhism asserts that this category is
equivalent to $\Rep(G^{\vee})$ as a tensor category. The corresponding fiber functor from $\Perv_{G(\calO)}(\Gr_G)$
to vector spaces sends every perverse sheaf $\calS\in\Perv_{G(\calO)}(\Gr_G)$ to its cohomology.

More precisely, one can show (we review it in \refs{grass}) that $G(\calO)$-orbits on $\Gr_G$ are finite-dimensional
and they are indexed by the set $\Lam^+$ of dominant weights of $G^{\vee}$. For every $\lam\in\Lam^+$ we denote
by $\Gr_G^{\lam}$ the corresponding orbit and by $\oGr_G^{\lam}$ its closure in $\Gr_G$. Then $\Gr_G^{\lam}$ is a
non-singular quasi-projective algebraic variety over $\CC$ and $\oGr_G^{\lam}$ is a (usually singular)
projective variety. One has
$$
\oGr_G^{\lam}=\bigcup\limits_{\mu\leq \lam}\Gr_G^{\lam}.
$$
One of the main properties of the geometric Satake isomorphism is that it sends the irreducible $G^{\vee}$-module
$L(\lam)$ to the intersection cohomology complex $\IC(\oGr_G^{\lam})$. In particular, the module $L(\lam)$
itself gets realized as the intersection cohomology of the variety $\oGr_G^{\lam}$.

As a byproduct of the geometric Satake isomorphism one can compute $\IC(\oGr_G^{\lam})$ in terms of the
module $L(\lam)$. Namely, it follows from \cite{Lu-qan}, \cite{Bry} and \cite{Gi} that the stalk
of $\IC(\oGr_G^{\lam})$ at a point of $\Gr_G^{\mu}$ as a graded vector space is essentially
equal to the associated graded $\gr^F L(\lam)_{\mu}$ of the $\mu$-weight space $L(\lam)_{\mu}$ in $L(\lam)$
with respect to certain filtration (called sometimes {\em the Brylinski-Kostant filtration}).
We refer the reader to \refss{q-an} for more details.

In fact, in \cite{KT} the authors construct certain canonical transversal slice $\ocalW^{\lam}_{G,\mu}$ to
$\Gr_G^{\mu}$ inside $\oGr_G^{\lam}$. This is a conical affine algebraic variety (i.e. it is endowed with an
action of the multiplicative group $\GG_m$ which contracts it to one point).
The above result about the stalks of $\IC(\oGr_G^{\lam})$
then gets translated into saying that the stalk of the IC-sheaf of $\ocalW^{\lam}_{G,\mu}$ at the unique $\GG_m$-fixex point
is essentially isomorphic to $\gr^F L(\lam)_{\mu}$. Note that $\ocalW^{\lam}_{G,\mu}$ is contracted to the above point by the
$\GG_m$-action, it follows that the stalk of of $\IC(\ocalW^{\lam}_{G,\mu})$ is equal to the global intersection cohomology
$\IH^*(\ocalW^{\lam}_{G,\mu})$.
\ssec{}{The affine case}This paper should be considered as a part of an attempt
to generalize Langlands duality to the case of affine Kac-Moody groups.
More precisely, in this paper we generalize some of the facts related to the geometric Satake isomorphism to the affine
case. Let us mention that simultaneously in \cite{BrKa} an attempt is made also to generalize the more classical story
(such as the usual Satake isomorphism) to the affine case.

Let us now assume that $G$ is semi-simple and simply connected and let $G_{\aff}$ denote the corresponding
(untwisted) affine Kac-Moody group (cf. \refs{snaive}); let $\grg_{\aff}$ be the Lie algebra of $G_{\aff}$.
In \refs{snaive} we explain that to $G_{\aff}$ one can attach
the corresponding dual affine Kac-Moody group $G_{\aff}^{\vee}$. The Lie algebra $\grg_{\aff}^{\vee}$
of $G_{\aff}^{\vee}$ is an affine
Kac-Moody algebra  whose Dynkin diagram is dual to that of $\grg_{\aff}$.
If $\grg$ is simply laced, the algebra $\grg_{\aff}^{\vee}$ is isomorphic to $\grg_{\aff}$. In general,
however, this might be a twisted affine Lie algebra (in particular, in general
$\grg_{\aff}^{\vee}\not\simeq (\grg^{\vee})_{\aff}$).

The idea of this paper came to us from a conversation with I.~Frenkel who suggested that integrable
representations of $G_{\aff}$ of level $k$ should be realized geometrically in terms of some moduli spaces
related to $G$-bundles on $\AA^2/\Gam_k$, where $\Gam_k$ is the group of roots of unity of order $k$ acting
on $\AA^2$ by $\zeta(x,y)=(\zeta x,\zeta^{-1}y)$.
\footnote{I.~Frenkel's suggestion is discussed also in \cite{Lic} whose relation to this paper is not
clear to us at the moment.}
In this paper we give some kind of rigorous meaning to this suggestion with a small modification: $G_{\aff}$ has
to be replaced by $G_{\aff}^{\vee}$.
\footnote{In fact, I.~Frenkel made this suggestion only when $G$ is simply laced, in which case there is almost no
difference between $G_{\aff}$ and $G_{\aff}^{\vee}$.}
Roughly speaking, this paper constitutes an attempt to make precise the following
principle. Let $\Bun_G(\AA^2)$ denote the moduli space of principal $G$-bundles on $\PP^2$ trivialized at the ``infinite"
line
$\PP^1_{\infty}\subset \PP^2$. This is an algebraic variety which has connected components parametrized by non-negative
integers, corresponding to different values of the second Chern class of the corresponding bundles.
Similarly, one can define $\Bun_G(\AA^2/\Gam_k)$ (cf. \refss{actionofgam} for the details). Very vaguely, the main idea
of this paper can be formulated in the following way:

\medskip
\noindent
{\bf The basic principle:}

1) The integrable representations of $G_{\aff}^{\vee}$ (not of $G_{\aff}$!) of level $k$
have to do with the geometry (e.g. intersection cohomology) of some varieties closely related to
$\Bun_G(\AA^2/\Gam_k)$.

2) This relation should be thought of as similar to the relation between finite-dimensional representations
of $G^{\vee}$ and the geometry of the affine Grassmannian $\Gr_G$.

\medskip
\noindent
We believe that 1) above has many different aspects. In this paper we concentrate on just one such aspect;
namely, in \refs{formulation} we explain how one can construct an analog of the varieties $\ocalW^{\lam}_{\mu}$
in the affine case (using the variety $\Bun_G(\AA^2/\Gam_k)$ as well as the corresponding {\em Uhlenbeck
compactification} of the moduli space of $G$-bundles - cf. \cite{bfg} and \refs{formulation}).
We conjecture that the stalks of IC-sheaves of these varieties are governed by the affine version of
$\gr^F L(\lam)_{\mu}$ (cf. \refco{Greenland}).

Let us note that in the case $G=SL(n)$ a slightly different approach
to 1) is discussed in \cite{Lic}. Also, I.~Grojnowski has informed us that he made some attempts to
construct $\Gr_{G_{\aff}}$ (apparently, using completely different methods).

We do not know how to prove our main conjecture in general; in the last
three sections of the paper we check our conjecture in several special cases. In particular, in \refs{large}
we prove all our conjectures in the limit $k\to \infty$ (cf. \refs{large} for the exact formulation).
In \refs{k=1} we prove a slightly weaker version of our main \refco{Greenland} in the case $k=1$;
the proof is based on the results of \cite{bfg}. Also in \refs{sln}
we prove again a slightly
weaker form of our conjecture for $G=\SL(N)$). Let us mention the main ingredient of that proof.
Let $\grg$ be a simply laced simple finite-dimensional Lie algebra. Then (by McKay correspondence)
to $\grg$ one can associate a finite subgroup $\Gam$ of $\SL(2,\CC)$.
Recall that H.~Nakajima (cf. e.g. \cite{nakajima})
gave a geometric construction of integrable $\grg_{\aff}$-modules of level $N$ using
certain moduli spaces which, roughly speaking, have to do
with vector bundles of rank $N$ on $\AA^2/\Gam$.
In particular, if $\grg=\fsl(k)$ it follows that

\smallskip
1) By H.~Nakajima the geometry of vector bundles of rank $n$ on $\AA^2/\Gam_k$ is related to integrable
modules over $\fsl(k)_{\aff}$ of level $N$.

2) By I.~Frenkel's suggestion (and by our \refco{Greenland}) the geometry of vector bundles
of rank $N$ on $\AA^2/\Gam_k$ is related to integrable
modules over $\fsl(N)_{\aff}$ of level $k$.

\smallskip
\noindent
On the other hand, in the representation theory of affine Lie algebras
there is a well-known relation, due to I.~Frenkel,
between integrable
modules over $\fsl(k)_{\aff}$ of level $N$ and integrable
modules over $\fsl(N)_{\aff}$ of level $k$. This connection is called {\em level-rank duality}; one of its aspects
is discussed in \cite{Fr}. It turns out that combining the results of \cite{Fr} with the results of
\cite{nakajima} one gets a proof of our main conjecture (in a slightly weaker form).

\medskip
\noindent
It is of course reasonable to ask why $G$-bundles on $\AA^2/\Gam_k$ have anything to do with the sought-for
affine Grassmannian of $G_{\aff}$. We don't have a satisfactory answer to this question, though some sort
of explanation (which would be too long to reproduce in the Introduction) is provided in \refs{snaive}.
Also we have been informed by E.~Witten that this phenomenon (at least when $G$ is simply laced)
has an explanation in terms of 6-dimensional conformal field theory (cf. \cite{Wit}).

In~\cite{BF} we explore other aspects of 1); in particular, we consider
(mostly conjectural) affine analogs of convolution of $G(\calO)$-equivariant
perverse sheaves on $\Gr_G$, the analog
of the so called Beilinson-Drinfeld Grassmannian, etc.
\ssec{Contents}{}This paper is organized as follows: in \refs{grass} we review the basic results about the affine Grassmannian
of a reductive group $G$ that we would like to generalize to the affine case.
In \refs{snaive} we discuss how the affine Grassmannian of $G_{\aff}$ looks like from the naive (i.e. set-theoretic and
not algebraic-geometric) point of view. In \refs{formulation} we describe our affine analogs of the transversal
slices $\ocalW^{\lam}_{G,\mu}$ and formulate the main conjectures  about them (linking them
to integrable representations of $G_{\aff}^{\vee}$). The next three sections are devoted to a verification of these conjectures
in a number of special cases. In particular, in \refs{sln} we prove most of them in the case $G=SL(N)$ and in \refs{large}
we prove all our conjectures when $k$ is large compared to the other parameters and in \refs{k=1} we prove almost all of our
conjectures for $k=1$.
\ssec{}{Acknowledgments}As was mentioned above, the general idea saying that $G$-bundles on
$\AA^2/\ZZ_k$ should have something
to do with integrable modules of level $k$ was explained to the first author by I.~Frenkel, to whom we express
our deepest gratitude. The authors would also like to thank
B.~Feigin, D.~Gaitsgory, D.~Kazhdan, H.~Nakajima and E.~Witten
for many very interesting discussions on the subject.
In addition our work was very much influenced by the paper \cite{MVy} by
I.~Mirkovi\'c and M.~Vybornov to whom we also want to express our
gratitude. A.~B.~ was partially supported by the NSF Grant DMS-0600851.
M.~F. is grateful to MSRI for the hospitality
and support; he was partially supported by the RFBR grant 09-01-00242 and the
Science Foundation of the SU-HSE award No.09-08-0008 and 09-09-0009.
\sec{grass}{Basic results about affine Grassmannian}
\ssec{21}{Definition}

Let $\calK=\CC((s))$, $\calO=\CC[[s]]$. By the {\it affine
Grassmannian} of $G$ we will mean the quotient
$\gg=G(\calK)/G(\calO)$. It is known (cf. \cite{BD,MV}) that
$\gg$ is the set of $\CC$-points of an ind-scheme over
$\CC$, which we will denote by the same symbol. Note that $\gg$ is defined
for any (not necessarily reductive) group $G$.

Let $\Lam=\Lam_G$ denote the coweight lattice of $G$ and let $\Lam^{\vee}$
denote the dual lattice (this is the weight lattice of $G$).
We let $2\rho_G^{\vee}$ denote the sum of the positive roots of $G$.

The group-scheme $G(\calO)$ acts on $\gg$ on the left and
its orbits can be described as follows.
One can identify the lattice $\Lam_G$ with
the quotient $T(\calK)/T(\calO)$. Fix $\lam\in\Lam_G$ and
let $s^{\lam}$ denote any lift of $\lam$ to $T(\calK)$.
Let $\gg^{\lam}$ denote the $\go$-orbit of $s^{\lam}$
(this is clearly independent of the choice of $\lam(s)$).
The following result is well-known:
\lem{gras-orbits}
\begin{enumerate}
\item
$$
\gg=\bigcup\limits_{\lam\in\Lam_G}\gg^{\lam}.
$$
\item
We have $\Gr_G^{\lam}=\Gr_G^{\mu}$ if an only if $\lam$ and $\mu$ belong
to the same $W$-orbit on $\Lam_G$ (here $W$ is the Weyl group of $G$). In particular,
$$
\gg=\bigsqcup\limits_{\lam\in\Lam^+_G}\gg^{\lam}.
$$
\item
For every $\lam\in\Lam^+$ the orbit
$\gg^{\lam}$ is finite-dimensional and its dimension is
equal to $\la\lam,2\rho_G^{\vee}\ra$.
\end{enumerate}
\elem
Let $\ogl$ denote the closure of $\gg^{\lam}$ in $\gg$;
this is an irreducible projective algebraic variety; one has
$\gg^{\mu}\subset \ogl$ if and only if $\lam-\mu$ is a sum of positive roots of
$G^{\vee}$.
We will denote by $\IC^{\lambda}$ the intersection
cohomology complex on $\ogl$. Let $\pgg$ denote the category of
$G(\calO)$-equivariant perverse sheaves on $\gg$. It is known
that every object of this category is a direct sum of the
$\IC^{\lam}$'s.

The group $\CC^*$ acts naturally on $G((s
))$ by ``loop rotation"; in other words,
any $a\in\CC^*$ acts on $g(s)\in G((s))$ by $g(s)\mapsto g(as)$. The fixed point
variety $(\gg^{\lam})^{\CC^*}$ is known to be isomorphic as a $G$-variety (via the action
of $G$ on $(\gg^{\lam})^{\CC^*}$ coming from the natural identification $G=G(\calO)^{\CC^*}$)
to $G/P_{\lam}$ where $P_{\lam}\subset G$ is the parabolic subgroup of G satisfying the following conditions:

\smallskip
1) $P_{\lam}$ contains the standard
Borel subgroup $B_+$;

2)  The Weyl group of the corresponding Levi subgroup $M_{\lam}$
is equal to the stabilizer $W_{\lam}$ of $\lam$ in $W$.

\smallskip
\noindent
In particular, $(\gg^{\lam})^{\CC^*}$
is a projective variety. Moreover, the $\CC^*$-action contracts all of $\gg^{\lam}$ to $(\gg^{\lam})^{\CC^*}$
and the corresponding map $p_\lam:\gg^{\lam}\to (\gg^{\lam})^{\CC^*}$ is a fiber bundle with fiber being isomorphic
to an affine space~\cite{BB}.

The affine Grassmannian $\Gr_G$ admits a canonical principal $G$-bundle $\eps:\tilGr_G\to\Gr_G$.
The set of $\CC$-points of $\tilGr_G$ is just $G(\calK)/G(\calO)_1$ where $G(\calO)_1$ is the kernel
of the natural (``evaluation at 0") map $G(\calO)\to G$.
This $G$-bundle is equivariant with respect to both $G(\calK)$ and the loop rotations.
The following lemma is easy and its proof is left to the reader:
\lem{grastorsor}
Let $x\in(\Gr_G)^{\CC^*}$ (here as before the super-script $\CC^*$ means ``fixed points of the loop
rotation"). Then $x\in (\Gr_G^{\lam})^{\CC^*}$ if and only if the homomorphism $\CC^*\to G$ obtained
via the action of $\CC^*$ in $\eps^{-1}(x)$ is conjugate to $\lam$.
\elem
\ssec{q-an}{IC-stalks and the $q$-analog of the weight multiplicity}
Since every $\gg^{\lam}$ is simply connected, it follows that $\IC^{\lam}|_{\gg^{\mu}}$ is a constant
complex, corresponding to some graded vector space $\tilV^{\lam}_\mu=\oplus (\tilV^{\lam}_{\mu})_i$.
We let $V^{\lam}_{\mu}=\tilV^{\lam}_{\mu}[-\la\mu,2\rho_{G}^{\vee}\ra]$. In other words,
$(V^{\lam}_{\mu})_i=(\tilV^{\lam}_{\mu})_{i-\la \mu, 2\rho_G^{\vee}\ra}$.

It follows from
\cite{Lu-qan} and \cite{Bry} (as well as from \cite{Gi} whose proof is somewhat more adapted to our
purposes) that the space $V^{\lam}_{\mu}$ can be described as follows.

Let $L(\lam)$ denote the finite-dimensional representation of $G^{\vee}$ with highest weight $\lam$.
Let
$$
\gd=\grn_+^{\vee}\oplus\grt^{\vee}\oplus\grn_-^{\vee}
$$
be the natural triangular decomposition. Fix a regular nilpotent element $e\in \grn_+^{\vee}$
and define an increasing filtration on $F^i L(\lam)$ by putting
\eq{filtration}
x\in F^i L(\lam)\Leftrightarrow e^i(x)=0.
\end{equation}

Thus we have the following
\th{ic-grassman}
1) The graded vector space $V^{\lam}_{\mu}$ lives only in even degrees.

2) We have a canonical isomorphism
$$
(V^{\lam}_{\mu})_{-2i}\simeq \gr^F _i L(\lam)_\mu
$$
(here $L(\lam)_\mu$ denotes the $\mu$ weight space in $L(\lam)$).
\eth
\ssec{aroundkost}{The $q$-analog of Kostant multiplicity formula}
 The second assertion of \reft{ic-grassman} can be rewritten as follows. Introduce the generating functions:
$$
\IC^{\lam}_{\mu}(q):=\sum\limits_{i}\dim(V^{\lam}_{\mu})_{-2i}q^i;\quad
^eC^{\lam}_{\mu}(q):=\sum\limits_{i} \dim\gr_i^F L(\lam)_{\mu} q^i.
$$
Then the second assertion of \reft{ic-grassman} says that
$\IC^{\lam}_{\mu}(q)={^eC^{\lam}_{\mu}}$. We shall refer to this function as
{\em the $q$-character} of $L(\lam)_{\mu}$.

The $q$-character of $L(\lam)_{\mu}$ can be given another purely combinatorial definition.
Namely, let $R$ denote the set of coroots of $G$ (= roots of $G^{\vee}$) and let $R^+$ be
the corresponding set of positive roots. Let $\Lam^{\pos}$ denote the subsemigroup of $\Lam$ consisting of elements of
the form $\sum n_j \alp_j$ where $n_j\in \ZZ_{\geq 0}$ and $\alp_j\in R^+$. For all $\beta\in \Lam^{\pos}$ let us
introduce the Kostant partition function of $K_{\beta}(q)$ of $\beta$ by:
$$
\sum_{\beta\in\Lam^{\pos}}K_\beta(q):=
\prod_{\alpha\in R^+}(1-qe^\alpha)^{-1},
$$
For each $\lam\in\Lam$ let us set $w\cdot\lam=w(\lam+\rho_G)-\rho_G$. Also let $\ell:W\to \ZZ_+$ denote the
length function.

Now we define
\eq{kostant}
C^\lambda_\mu(q):=
\sum_{w\in W}(-1)^{\ell(w)}K_{w\cdot\lambda-\mu}(q).
\end{equation}

Then we have the following
\th{lusztig-brylinski}For every $\lam,\mu\in\Lam^+$
We have $\IC^{\lam}_{\mu}(q)={^eC^{\lam}_{\mu}(q)}=C^{\lam}_{\mu}(q)$.
\eth
The equality $\IC^{\lam}_{\mu}(q)=C^{\lam}_{\mu}(q)$ has been proved in Lusztig's paper
\cite{Lu-qan} where the author also formulated the equality ${^eC^{\lam}_{\mu}(q)}=C^{\lam}_{\mu}(q)$
as a conjecture. Note that this conjecture is purely representation-theoretic (i.e. it has nothing to do
with the geometry of $\Gr_G$); it was proved later by  Brylinski \cite{Bry}.
Note however that the equality $\IC^{\lam}_{\mu}(q)={^eC^{\lam}_{\mu}(q)}$ has an independent proof
(cf. \cite{Gi} in which an isomorphism as in \reft{ic-grassman} is constructed). Note also that both
${^eC^{\lam}_{\mu}(q)}$ and $C^{\lam}_{\mu}(q)$ make sense when $\mu$ is arbitrary (i.e. not necessarily).
In this case it was shown in \cite{joseph} that the equality ${^eC^{\lam}_{\mu}(q)}=C^{\lam}_{\mu}(q)$
still holds up to certain power of $q$ which depends on how
``non-dominant" $\mu$ is (cf.
Theorem 7.6 in \cite{joseph}).
\ssec{trans-finite}{Transversal slices}Consider the group $G[s^{-1}]\subset G((s))$; let us denote by
$G[s^{-1}]_1$ the kernel of the natural (``evaluation at $\infty$") homomorphism
$G[s^{-1}]\to G$. For any $\lam\in\Lam$ let $\Gr_{G,\lam}=G[s^{-1}]\cdot s^{\lam}$. Then
it is easy to see that one has
$$
\gg=\bigsqcup\limits_{\lam\in\Lam^+}\Gr_{G,\lam}
$$

Let also $\calW_{G,\lam}$ denote the $G[s^{-1}]_1$-orbit of $s^{\lam}$.
For any $\lam,\mu\in\Lam^+$, $\lam\geq \mu$ set
$$
\Gr^{\lam}_{G,\mu}=\gg^{\lam}\cap \Gr_{G,\mu},\quad
\oGr^{\lam}_{G,\mu}=\oGr_G^{\lam}\cap \Gr_{G,\mu}
$$
and
$$
\calW^{\lam}_{G,\mu}=\gg^{\lam}\cap \calW_{G,\mu},\quad
\ocalW^{\lam}_{G,\mu}=\oGr_G^{\lam}\cap \calW_{G,\mu}.
$$
Note that $\ocalW^{\lam}_{G,\mu}$ contains the point $s^{\mu}$ in it.
\lem{}
\begin{enumerate}
\item
The point $s^{\mu}$ is the only $\CC^*$-fixed point in $\ocalW^{\lam}_{G,\mu}$.
The action of $\CC^*$ on $\ocalW^{\lam}_{G,\mu}$ is ``repelling", i.e. for any
$w\in \ocalW^{\lam}_{G,\mu}$ we have
$$
\lim\limits_{a\to\infty} a(w)=s^{\mu}.
$$
\item
The variety $\oGr^{\lam}_{G,\mu}$ is a fiber bundle over $G/P_{\mu}$ with fiber
$\ocalW^{\lam}_{G,\mu}$.
\item
There exists an open subset $U$ in $\Gr_G^{\mu}$ and an open embedding
$U\x \ocalW^{\lam}_{G,\mu}\hookrightarrow \oGr^{\lam}_G$ such that the diagram
$$
\begin{CD}
U\x \{s^{\mu}\} @>>> \Gr^{\mu}_{G}\x\{s^{\mu}\}\\
@VVV @VVV\\
U\x \ocalW^{\lam}_{G,\mu} @>>> \oGr^{\lam}_G
 \end{CD}
$$
is commutative.
In other words, $\ocalW^{\lam}_{G,\mu}$ is a transversal slice to $\Gr^{\mu}_G$ inside
$\oGr^{\lam}_G$. In particular, the stalk of $\IC_{\ocalW^{\lam}_{G,\mu}}$
at the point $s^{\mu}$ is
equal to $V^{\lam}_{\mu}[-\la \mu, 2\rho^{\vee}\ra]$ (note that $\la \mu, 2\rho^{\vee}\ra$ is the
dimension on $\Gr^{\mu}_G$).
\end{enumerate}
\elem

\prf
The first two statements are obvious, and the third one follows from
Propositions 1.3.1 and 1.3.2 of \cite{KT}.
\epr

\sec{snaive}{The double affine Grassmannian: the naive approach}
\ssec{}{The affine Kac-Moody group and its Langlands dual}
In this paper for convenience we adopt the polynomial version of loop groups (as opposed to formal
loops version).

Let $G$ be a connected simply connected algebraic group with simple Lie algebra $\grg$.
Consider the corresponding polynomial algebra $\grg[t,t^{-1}]$ and the group $G[t,t^{-1}]$.
It is well-known that $G[t,t^{-1}]$ has a canonical central extension $\hatG$:
$$
1\to \CC^*\to \hatG\to G[t,t^{-1}]\to 1.
$$
We denote by $\hgrg$ the corresponding Lie algebra; it fits into the exact sequence
$$
0\to \CC\to \hgrg\to\grg[t,t^{-1}]\to 0.
$$

This central extension corresponds to a choice of an invariant symmetric bilinear form
$(\cdot,\cdot)$ on $\grg$ which is integral and even on $\Lam\subset \grt\subset \grg$.
We shall assume that $(\cdot,\cdot)$ is the minimal possible such form (i.e.
that $(\alp,\alp)=2$ for any short coroot $\alp$ of $G$).

The group $\CC^*$ acts naturally on $G[t,t^{-1}]$ and this action lifts to $\hatG$.
We denote the corresponding semi-direct product by $G_{\aff}$; we also let $\grg_{\aff}$ denote its Lie algebra.
Thus $\grg_{\aff}$ is an untwisted affine Kac-Moody Lie algebra,
in particular, it can be described by the corresponding affine root system.
We let $G'_{\aff}$  be the quotient of $G_{\aff}$ by the central $\CC^*$.

Let now $\grg_{\aff}^{\vee}$ denote the {\it Langlands dual} affine Lie algebra. By definition, this is an affine
Kac-Moody Lie algebra whose root system is dual to that of $\grg_{\aff}$. Note that in general (when $\grg$ is not
simply laced) the algebra $\grg_{\aff}^{\vee}$ is not isomorphic to $(\grg^{\vee})_{\aff}$ (here $\grg^{\vee}$ as before denotes
the Langlands dual Lie algebra of $\grg$). Moreover, if $\grg$ is not simply laced, then $\grg_{\aff}^{\vee}$ is a twisted
affine Lie algebra. However, the algebra $\grg_{\aff}^{\vee}$ always contains $\grg^{\vee}\x \CC^2$ as a Levi subalgebra.

Let also $G_{\aff}^{\vee}$ be the corresponding Langlands dual group. This is a connected
affine Kac-Moody group in the sense
of \cite{Kum} which is uniquely characterized by the following two properties:

$\bullet$ The Lie algebra of $G_{\aff}^{\vee}$ is $\grg_{\aff}^{\vee}$;

$\bullet$ The embedding $\grg^{\vee}\hookrightarrow \grg_{\aff}^{\vee}$ lifts to an embedding
$G^{\vee}\hookrightarrow G_{\aff}^{\vee}$.

The group $G_{\aff}^{\vee}$ maps naturally to $\CC^*$ (this homomorphism is dual to the central embedding
$\CC^*\to G_{\aff}$). We denote the kernel of this homomorphism by $\hatG^{\vee}$ and we let
$\hgrg^{\vee}$ denote its Lie algebra.
\ssec{vee}{The weight lattice of $G_{\aff}^{\vee}$}
In what follows we shall write just $\Lam$ instead of $\Lam_G$; also from now on we shall usually denote elements
of $\Lam$ by $\blambda,\bmu...$ (instead of just writing $\lam,\mu...$ in order to distinguish them from
the coweights of $G_{\aff}$ (= weights of $G_{\aff}^{\vee}$).

The weight lattice of $G_{\aff}^{\vee}$ (resp. of
$\hatG^{\vee}$) is naturally identified with $\Lam_{\aff}=\ZZ\x\Lam\x \ZZ$ (resp. with $\hatLam=\ZZ\x\Lam$).
Here the first $\ZZ$-factor
is responsible for the center of $G_{\aff}^{\vee}$ (or $\hatG^{\vee}$);
it can also be thought of as coming from the loop
rotation in $G_{\aff}$. The second $\ZZ$-factor is responsible for the loop rotation in $G_{\aff}^{\vee}$ (and thus it is absent
in the case of $\hatG^{\vee}$; it may also be thought of
as coming from the center of $G_{\aff}$).
We denote by $\Lam_{\aff}^+$ the set of dominant weights of $G_{\aff}^{\vee}$ (which is the same as the set of dominant
coweights of $G_{\aff}$). We also denote by $\Lam_{\aff,k}$ the set of weights of $G_{\aff}^{\vee}$ of level $k$,
i.e. all the weights of the form $(k,\olam,n)$. We put $\Lam_{\aff,k}^+=\Lam_{\aff}^+\cap \Lam_{\aff,k}$.

Let
$\Lam_k^+\subset \Lam$ denote the set of dominant coweights of $G$ such that $\la \blambda,\alp)\leq k$
when $\alp$ is the highest root of $\grg$.
Then it is well-known that a weight $(k,\olam,n)$ of
$G_{\aff}^{\vee}$ lies in $\Lam_{\aff,k}^+$  if and only if $\olam\in\Lam_k^+$ (thus $\Lam_{\aff,k}=\Lam_k^+\x \ZZ$).

Let also $W_{\aff}$ denote affine Weyl group of $G$ which is the semi-direct product of $W$ and $\Lam$.
It acts on the lattice $\Lam_{\aff}$ (resp. $\hatLam$)
preserving each $\Lam_{\aff,k}$ (resp. each $\hatLam_k$). In order to describe this action explicitly it is convenient
to set $W_{\aff,k}=W\ltimes k\Lam$ which naturally acts on $\Lam$. Of course the groups $W_{\aff,k}$ are canonically
isomorphic to $W_{\aff}$ for all $k$. Then the restriction of the $W_{\aff}$-action to $\Lam_{\aff,k}\simeq\Lam\x\ZZ$
comes from the natural $W_{\aff,k}$-action on the first multiple.

It is well known that every $W_{\aff}$-orbit on $\Lam_{\aff,k}$ contains unique element of $\Lam_{\aff,k}^+$.
This is equivalent to saying that $\Lam_k^+\simeq \Lam/W_{\aff,k}$.

It is clear that for two integers $0\leq k<l$ we have $\Lam_k^+\subset
\Lam_l^+$. However, if $k|l$, we also have the natural
embedding $W_{\aff,l}\subset W_{\aff,k}$ and thus the identification $\Lam_k^+=\Lam/W_{\aff,k}, \Lam_l^+=\Lam/W_{\aff,l}$
gives rise to the projection $\Lam_l^+\to \Lam_k^+$ which is equal to identity on $\Lam_k^+$.

Let $\Gam_k$ denote the group of roots of unity of order $k$. Then we claim the following:
\lem{}There is a natural bijection between $\Lam/W_{\aff,k}$ and the set of conjugacy classes of
homomorphisms $\Gam_k\to G$.
\footnote{This lemma is true for an arbitrary reductive group $G$ (with the same proof). However, the identification
$\Lam_k^+=\Lam/W_{\aff,k}$ is only true when $G$ is simply connected.}
\elem
\prf
Recall that $\Lam=\Hom(\GG_m,T)$. Thus, since $\Gam_k$ is a closed subscheme of $\GG_m$,
given every $\olam\in \Lam$ we may restrict it to $\Gam_k$ and get a homomorphism $\Gam_k\to T$.
By composing it with the embedding $T\hookrightarrow G$ we get a homomorphism $\Gam_k\to G$ which clearly
depends only on the image of $\lam$ in $\Lam/W_{\aff,k}$. Thus we get a well-defined map
$\Lam/W_{\aff,k}\to (\Hom(\Gam_k,G)/G)$. The surjectivity of this map follows from the fact that
$\Gam_k$ is finite and thus diagonalizable. For the injectivity, note that any two elements in $T$ which are
conjugate in $G$ lie in the same $W$-orbit in $T$. Thus it is enough to show that
for any two homomorphisms $\lam,\mu:\GG_m\to T$ whose restrictions to $\Gam_k$ coincide the difference
$\lam-\mu$ is divisible by $k$. This is enough to check for $T=\GG_m$ where it is obvious.
\epr
\ssec{referee}{Representations of $G_{\aff}^{\vee}$}
In the sequel we shall be concerned with representations of $G_{\aff}^{\vee}$ (i.e. with integrable
representations of the Lie algebra $\grg_{\aff}^{\vee}$ which integrate to $G_{\aff}^{\vee}$). The following
results are well-known (cf. \cite{Kum}):

a) The irreducible representations of $\hatG^{\vee}$ of level $k$ are in one-to-one correspondence with elements of $\Lam_k^+$;

b) The irreducible representations of $G_{\aff}^{\vee}$ of level $k$
are in one-to-one correspondence with elements of $\ZZ\x\Lam_k^+$;

c) For any irreducible representation $L$ of $G_{\aff}^{\vee}$, its restriction to $\hatG^{\vee}$ is irreducible.

For every $\lam\in\Lam_{\aff,k}^+$ we denote by $L(\lam)$ the corresponding representation of $G_{\aff}^{\vee}$.
For every $\mu\in\Lam_{k,\aff}$ we denote by $L(\lam)_{\mu}$ the $\mu$-weight space of $L(\lam)$.
It is known to be finite-dimensional.

Recall that the Lie algebra $\grg_{\aff}^{\vee}$ has a triangular decomposition
$$
\grg_{\aff}^{\vee}=\grn_{-,\aff}^{\vee}\oplus\grt_{\aff}^{\vee}\oplus\grn_{+,\aff}^{\vee}.
$$
The Lie algebra $\grn_{+,\aff}^{\vee}$ is generated by the standard Cartan generators $e_0,e_1,...,e_r$.
Set $e=e_0+...+e_r$. This is
an affine analogue of a regular nilpotent element considered in \refs{grass}. We define a filtration $F^iL(\lam)$ in
the same way as in \refe{filtration}; it can be restricted to every
weight space
$L(\lambda)_{\mu}$.
We set $\gr^F L(\lambda)_\mu$ to be the associated graded space.
\ssec{}{The affine Grassmannian of $G_{\aff}$: the naive approach}The rest of this Section is written for motivational
purposes only; it will not be formally used in the rest of the paper.

Let $\calK=\CC[s,s^{-1}]$, $\calO=\CC[s]$ (note that this notation is different from the one we used in~\refss{21}).

Then we may consider the group $G_{\aff}(\calK)$
and its subgroup $G_{\aff}(\calO)$.

Let $\pi$ denote the natural map from $G_{\aff}(\calK)$ to $\calK^*=\CC^*\x \ZZ$.
We denote by $G_{\aff}^+(\calK)$ the sub-semigroup of $G_{\aff}(\calK)$ defined as follows:
$$
G_{\aff}^+(\calK)=G_{\aff}(\calO)\cup\{ g\in G_{\aff}(\calK)\ \text{such that
$\pi(g)=as^k$ where $k\geq0$}\}.
$$
For each $k\geq 0$ we consider the subset $G_{\aff,k}(\calK)$ of $G_{\aff}^+$ consisting of those $g\in G_{\aff}^+$
for which $\pi(g)=as^k$ (note that for all $k\neq 0$ we may replace $G_{\aff}^+$ by $G_{\aff}$ here).

The same definitions make sense when $G_{\aff}$ is replaced by $G'_{\aff}$.
\lem{naive}
The set $G_{\aff}(\calO)\backslash G_{\aff,k}(\calK)/G_{\aff}(\calO)$ is in natural bijection with $\Lam_{\aff,k}^+$. Similarly,
$G'_{\aff,k}(\calO)\backslash G_{\aff,k}'(\calK)/G_{\aff}'(\calO)$ is in natural bijection with $\Lam_k^+$.
\footnote{If one uses $\Lam_{\aff,k}/W_{\aff}$ and $\Lam/W_{\aff,k}$ instead $\Lam_{\aff,k}^+$ and $\Lam_k^+$ the statement
becomes true for arbitrary reductive $G$.}
\elem
Let us formulate \refl{naive} a bit more precisely.
Recall that $G_{\aff}$ contains $\GG_m\x G\x \GG_m$ as a subgroup (here the first $\GG_m$-factor stands for
the center and the second one --- for the loop rotation). In particular, $G_{\aff}(\calK)$ contains
$\calK^*\x G(\calK)\x \calK^*$. For every $\lam\in\Lam_{\aff,k}$ of the form $(k,\olam,n)$ we may consider the
element $s^{\lam}=(s^n, s^{\olam},s^k)$ in $\calK^*\x T(\calK)\x \calK^*\subset G_{\aff,k}(\calK)$.

We claim that
\begin{enumerate}
\item
$G_{\aff,k}(\calK)=\bigcup\limits_{\lam\in\Lam_{\aff,k}}G_{\aff}(\calO)\cdot s^{\lam}\cdot G_{\aff}(\calO)$;
\item
We have $G_{\aff}(\calO)\cdot s^{\lam}\cdot G_{\aff}(\calO)=G_{\aff}(\calO)\cdot
s^{\lam'}\cdot G_{\aff}(\calO)$ if and only if
$\lam$ and $\lam'$ are in the same orbit of $W_{\aff}$.
\end{enumerate}
Similar statements hold for $G_{\aff}'$ instead of $G_{\aff}$.

This is a direct ``affine" analog of the first two assertions of \refl{gras-orbits} and for $k=0$ it follows from there.
We shall be mostly concerned with the case $k>0$. In this case one can prove the above two assertions
directly; however, for the future it will be more instructive to prove it using the following trick. We shall do it
for $G_{\aff}'$ instead of $G_{\aff}$ (it is easy to see that our statements for $G_{\aff}$ and $G_{\aff}'$ are equivalent).

Let $\oS_k=\Spec \CC[x,y,z]/xy-z^k$ and let $S_k$ be the complement to the point
$(0,0,0)$ in $\oS_k$. This is a smooth surface.

The natural map
$p_k:\AA^2\to \oS_k$ given by
$$
(u,v)\mapsto (u^k,v^k,uv)
$$
identifies
$S_k$ with $(\AA^2\backslash \{ 0\})/\Gam_k$, where $\Gam_k$ acts on $\AA^2$ by
$\zeta(u,v)=(\zeta v,\zeta^{-1}v)$
(note that the action of $\Gam_k$ on $(\AA^2\backslash \{ 0\})$ is free).

For any $\CC$-variety $S$ we shall denote by $\Bun_G(S)$ the set of isomorphism classes of principal $G$-bundles
on $S$.

\prop{bijectionw}
The following sets are in natural bijection:
\begin{enumerate}
\item
The set $\Bun_G(S_k)$;
\item
The set $G'_{\aff}(\calO)\backslash G'_{\aff,k}/G'_{\aff}(\calO)$;
\item
The set of $G$-conjugacy classes of homomorphisms $\Gam_k\to G$.
\end{enumerate}
\eprop
\prf
Let us first establish the bijection between (1) and (3). Let $\calF$ denote a $G$-bundle on $S_k$.
Consider the $G$-bundle $p_k^*(\calF)$ on $\AA^2\backslash\{ 0\}$. It extends uniquely to the whole
of $\AA^2$ and thus it is trivial. On the other hand, the bundle $p_k^*(\calF)$ is $\Gam_k$-equivariant.
Since this bundle is trivial, such an equivariant structure gives rise to a homomorphism $\Gam_k\to G$
defined uniquely up to conjugacy. This defines a map $\Bun_G(S_k)\to \Hom(\Gam_k,G)/G$. It is clear that
this is actually a bijection, since a $\Gam_k$-equivariant $G$-bundle on $\AA^2\backslash\{ 0\}$
descends uniquely to $S_k$.

Let us now establish the bijection between (1) and (2). To do that let us denote by $\sig$ the automorphism
of $G(\calK)$ sending $g(t,s)$ to $g(ts,s)$. Let us now identify $\pi^{-1}(s^k)$ with $G(\calK)$ with right
$G(\calK)$ action being the standard one (by right shifts) and with left $G(\calK)$-action given by
$g(h)=\sig^k(g)h$. Thus we have the bijection
\eq{maindoubleq}
G'_{\aff}(\calO)\backslash G'_{\aff,k}/G'_{\aff}(\calO)=
G[ts^k,t^{-1}s^{-k},s]\backslash G[t,t^{-1},s,s^{-1}]/G[t,t^{-1},s].
\end{equation}
Let $U_k=\Spec \CC[ts^k,t^{-1}s^{-k},s]$, $V_k=\Spec \CC[t,t^{-1},s]$. Both $U_k$ and $V_k$ are isomorphic to
$\GG_m\x \AA^1$ and thus every $G$-bundle on either of these surfaces is trivial. Both $U_k$ and $V_k$ contain
$W=\Spec \CC[t,t^{-1},s,s^{-1}]\simeq \GG_m\x \GG_m$ as a Zariski open subset.
Thus the RHS of \refe{maindoubleq} can be identified with $\Bun_G(S_k')$ where $S_k'$ is obtained by gluing
$U_k$ and $V_k$ along $W$. Hence it remains to construct an isomorphism
$S_k'{\widetilde \to} S_k$. Such an isomorphism can be obtained by setting $x=ts^k, y=t^{-1}, z=s$.
\epr

The affine Grassmannian $G_{\aff}(\calK)/G_{\aff}(\calO)$ of $G_{\aff}$ is a $\ZZ$-torsor over the corresponding
quotient for $G'_{\aff}$.
It is not clear how to think about the affine Grassmannian of $G_{\aff}$ in terms of algebraic geometry.
Let us note though that for each $k>0$ the set of double cosets
$G_{\aff}(\calO)\backslash G_{\aff,k}(\calK)/G_{\aff}(\calO)$
is also a $\ZZ$-torsor over the set $G'_{\aff}(\calO)\backslash G'_{\aff,k}(\calK)/G'_{\aff}(\calO)=\Lam^+_k$.
Moreover, by analyzing the geometric construction of this torsor given
in~\cite{BrKa} more carefully
and using the embedding $\AA^2\subset \PP^2$ one can show it follows from the above that the above
$\ZZ$-torsor can be canonically trivialized. Thus, if we set
$\Gr_{G_{\aff},k}=G_{\aff,k}(\calK)/G_{\aff}(\calO)$, we see that $G_{\aff}(\calO)$-orbits on
$\Gr_{G_{\aff},k}$ are parameterized by $\Lam_k^+\x\ZZ=\Lam_{\aff,k}^+$.
Alternatively, $G_{\aff}(\calO)$-orbits on
$\Gr_{G_{\aff},k}$ are parameterized by $\Bun_G(S_k)\x\ZZ$.
This remark might serve as a motivation for the constructions of next Section.

\sec{formulation}{The main construction}
\ssec{}{The moduli space of $G$-bundles on $\AA^2$}
Here we follow the exposition of \cite{bfg}. Let $\bfS$ be a smooth projective surface
containing $\AA^2$ as an open subset and let $\bfD_{\infty}$ denote the complement of
$\AA^2$ in $\bfS$ (we shall refer to it as ``the divisor at $\infty$"). In what follows
$\bfS$ will almost always be $\PP^2$, or sometimes $\PP^1\x \PP^1$.
For an integer $a$ let $\Bun_G^a(\AA^2)$ denote the moduli space of principal $G$-bundles on
$\bfS$ of second Chern class $a$ with a chosen trivialization on $\bfD_{\infty}$
 (cf. \cite{bfg} for the discussion of the notion of $c_2$ in this case). It is easy to see that
 $\Bun_G^d$ does not depend on the choice of $\bfS$. It is shown in \cite{bfg} that
 this space has the following properties:

 a) $\Bun_G^a(\AA^2)$ is non-empty if and only if $a\geq 0$;

 b) For $a\geq 0$ the space $\Bun_G^a(\AA^2)$ is an irreducible smooth
 quasi-affine variety of dimension $2a\check{h}$ where $\check{h}$ denotes the
dual Coxeter number of $G$.
\ssec{}{The Uhlenbeck space of $\AA^2$}
In \cite{bfg} we construct an affine scheme $\calU^a_G(\AA^2)$ containing $\Bun_G^a(\AA^2)$
as a dense open subset which we are going to  call the {\it Uhlenbeck space} of bundles on $\AA^2$.
\footnote{In fact in \cite{bfg} we give several definitions of $\calU^a_G(\AA^2)$ about which we only know that
the corresponding reduced schemes coincide. In this paper we shall take $\calU^a_G(\AA^2)$ to be reduced by definition.}

The scheme $\calU^a_G(\AA^2)$ is still irreducible but in general it is highly singular. The main property of $\calU^a_G(\AA^2)$
is that it possesses the following stratification:
\eq{stratification}
\Bun_G^a(\AA^2)=\bigcup\limits_{0\leq b\leq a}\Bun_G^b(\AA^2)\x \Sym^{a-b}(\AA^2).
\end{equation}
Here each $\Bun_G^b(\AA^2)\x \Sym^{a-b}(\AA^2)$ is a locally closed subset of $\calU^a_G(\AA^2)$ and its closure is
equal to the union of similar subsets corresponding to all $b'\leq b$.

We shall denote by $\Bun_G(\AA^2)$ (resp. $\calU_G(\AA^2)$) the disjoint union of all $\Bun_G^a(\AA^2)$
(resp. of $\calU^a_G(\AA^2)$).

Let us note that the group $G\x \GL(2)$ acts naturally on $\Bun_G^a(\AA^2)$: here the first factor acts by changing
the trivialization on $\bfD_{\infty}$ and the second factor acts on $\AA^2$ (formally we should take $\bfS=\PP^2$ for this;
then $\GL(2)$ acts on $\bfS$ preserving $\bfD_{\infty}$ and thus acts on $\Bun_G^a(\AA^2)$).
It is easy to deduce from the construction of \cite{bfg} that this action extends to an action of the same group
on the Uhlenbeck space $\calU^a_G(\AA^2)$.

The following remark (also proved in \cite{bfg}) will be needed in the future: consider the central $\CC^*\subset \GL(2)$.
Then its acts on $\calU^a_G(\AA^2)$ with unique fixed point $0^a\in\calU^a_G(\AA^2)$. In terms of the stratification
\refe{stratification} the point $0^a$ corresponds to the case $b=0$ (recall that $\Bun^0_G(\AA^2)$ consists of just one
point (the trivial bundle)) and to the point $(0,0)$ taken with multiplicity $a$ (as a point in
$\Sym^a(\AA^2)$. Moreover, the above $\CC^*$ action contracts all of $\calU^a_G(\AA^2)$ to $0^a$ (i.e.
for any $p\in\calU^a_G(\AA^2)$ we have $\lim_{z\to 0} z(p)=0^a$. In particular, any closed $\CC^*$-invariant
subvariety $W$ of $\calU^a_G(\AA^2)$ must also contain the point $0^a$ and all of $W$ is contracted to
$0^a$ by the $\CC^*$-action.
\ssec{factor}{The factorization morphism}
Let $\AA^{(a)}$ denote the $a$-th symmetric power of $\AA^1$. This is the space of all effective divisors
on $\AA^1$ of degree $a$. As a scheme, it is, of course, isomorphic to the affine space $\AA^a$.

In \cite{bfg} we construct the {\it factorization morphism}
$$
\pi^a:\Bun^a_G(\AA^2)\to \AA^{(a)}.
$$
This morphism, in fact, depends on the choice of a direction in $\AA^2$. Below we are going to work with
the ``horizontal" factorization morphism (in the terminology of \cite{bfg}), i.e. our $\AA^{(a)}$ should be thought
of as the $a$-th symmetric power of the horizontal axis in $\AA^2$. To describe this morphism, it is convenient
to take $\bfS=\PP^1\x \PP^1$. Let $\calF\in\Bun_G^a(\AA^2)$. Then $x\in\AA^1$ lies in the support of the divisor
$\pi^a(\calF)$ if and only if the restriction of $\calF$ to the line $\{ a\}\x\PP^1$ is non-trivial (i.e. $\{ a\}\x\PP^1$ is a ``jumping line'' for $\calF$).
We refer the reader to \cite{bfg} for the formal definition of $\pi^a$.

The term ``factorization morphism" comes from the following property of $\pi^a$:
assume that we are given two effective divisors
$D_1$ and $D_2$ on $\AA^1$ of degrees $a_1$ and $a_2$. Assume that $D_1$ and $D_2$ are
disjoint; in other words let assume that the supports of $D_1$ and $D_2$ do not intersect.
Let $a=a_1+a_2$. Then $\pi^a$ enjoys the following
{\it factorization property}:
\eq{factor}
\text{There a canonical isomorphism}\ (\pi^a)^{-1}(D_1+D_2)=(\pi^{a_1})^{-1}(D_1)\x (\pi^{a_2})^{-1}(D_2).
\end{equation}
The factorization morphism $\pi^a$ extends to $\calU^a_G(\AA^2)$; abusing the notation we shall denote it by
the same symbol. The factorization property \refe{factor} holds for this new $\pi^a$ as well.
\ssec{actionofgam}{Action of $\Gam_k$}
Let $\Gam_k$ as before denote the group of roots of unity of order $k$. We consider the embedding
$\Gam_k\hookrightarrow(\CC^*)^2\subset \GL(2)$ sending every $\zeta\in\Gam_k$ to $(\zeta,\zeta^{-1})$.
We shall refer to the corresponding action of $\Gam_k$ on $\AA^2$ as the {\it symplectic action} (since
it preserves the natural symplectic form on $\AA^2$). In what follows we shall always take $\bfS$ to be either
$\PP^2$ of $\PP^1\x \PP^1$; in both cases the symplectic action of $\Gam_k$ on $\AA^2$ extends to an action on
$\bfS$.

Since the group $G$ acts on $\Bun_G(\AA^2)$ (cf. previous subsection) it follows that a
choice of $\bmu\in\Lam^+_k$ (which as before we shall identify with the set of conjugacy classes
of homomorphisms $\Gam_k\to G$) gives rise to a homomorphism $\Gam_k\to G\x\GL(2)$ and thus to an
action of $\Gamma_k$ on
$\Bun_G(\AA^2)$ (strictly speaking for this we have to choose a lifting of $\bmu$ to $\Hom(\Gam_k,G)$;
however, for different choices of this lifting (lying automatically in the same conjugacy class)
we obtain isomorphic actions of $\Gam_k$).
We denote the fixed point set of this action by
$\Bun_{G,\bmu}(\AA^2/\Gamma_k)$.

It follows from the fact the elements of $\Bun_G(\AA^2)$ have no non-trivial automorphisms that
every $\calF\in \Bun_{G,\bmu}(\AA^2/\Gamma_k)$ is actually $\Gam_k$-equivariant bundle on
$\bfS=\PP^2$.
In particular, for every
$\CF\in \Bun_{G,\bmu}(\AA^2/\Gamma_k)$ the group $\Gam_k$  acts on the fiber
of $\CF_0$ of $\calF$ at the $\Gamma_k$-fixed point $0\in\AA^2$. Hence $\calF_0$
defines a conjugacy class of maps $\Gam_k\to G$, i.e. an element of $\Lam^+_k$.
For every $\blambda\in\Lam^+_k$ we denote by
$\Bun_{G,\bmu}^{\blambda}(\AA^2/\Gamma_k)$ the subset of $\Bun_{G,\bmu}(\AA^2/\Gam_k)$
formed by all
$\CF\in \Bun_{\bmu}(\AA^2/\Gamma_k)$ such that $\CF_0$ if of class $\blambda$.
Clearly, it is
a union of connected components of the fixed point set
$\Bun_{\bmu}(\AA^2/\Gamma_k)$.
We denote by $\Bun_{G,\bmu}^{\blambda,a}(\AA^2/\Gamma_k)$ the intersection of
$\Bun_{G,\bmu}^{\blambda}(\AA^2/\Gamma_k)$ with $\Bun_G^a(\AA^2)$.

\conj{connect}
$\Bun_{G,\bmu}^{\blambda,a}(\AA^2/\Gamma_k)$ is connected (possibly empty).
\econj
\ssec{}{Remark}
Most of the results below rest upon this conjecture, so we assume it
in what follows. In fact, later we are going to see that \refco{connect} holds in the following cases:

1) In the case $k=1$ (by \cite{bfg}).

2)  In the case $G=\SL(N)$. In this case \refco{connect} follows from the well-known results of W.~Crawley-Boevey
on connectedness of some quiver varieties~\cite{CB}, cf.~\refs{sln}
for the details.

3) For $k\ggg0$ (i.e. if we fix $a$, $\olam$ and $\omu$ and make
$k$ sufficiently large).

4) In the case when $\olam\geq \omu$ we can't quite prove \refco{connect} in full generality but we shall show in
\refs{large} that in this case we can single out some specific connected component of
$\Bun_{G,\bmu}^{\blambda,a}(\AA^2)$ which is, in fact,
sufficient for our purposes.
\defe{uhl}
\begin{enumerate}
\item
We define $\calU_{G,\bmu}^{\blambda,a}(\AA^2/\Gamma_k)$ as the closure of
$\Bun_{G,\bmu}^{\blambda,a}(\AA^2/\Gamma_k)$ inside $\calU_G^a(\AA^2)$.
\item
Let $\mu=(k,\bmu,m),\ \lambda=(k,\blambda,l)$ be two elements of $\Lam_{\aff}^+$. Then we set
\eq{deflammu}
\Bun^{\lambda}_{G,\mu}(\AA^2/\Gamma_k)=\Bun^{\blambda,a}_{G,\bmu}(\AA^2/\Gamma_k)\
\text{where $a=k(l-m)+\frac{(\blambda,\blambda)}{2}-\frac{(\bmu,\bmu)}{2}$}.
\end{equation}
Similarly, we define $\calU^{\lam}_{G,\mu}(\AA^2/\Gam_k)$ as the closure
$\Bun_{G,\mu}^{\lambda}(\AA^2/\Gamma_k)$ inside $\calU_G^a(\AA^2)$.
\end{enumerate}
\edefe
The choice of $a$ in \refe{deflammu} might seem a bit bizarre at the first glance; it can be partly justified
by the following result:
\th{dimension}
The dimension of $\Bun^{\lam}_{G,\mu}(\AA^2/\Gam_k)$ (and thus of $\calU_{G,\mu}^{\lam}(\AA^2/\Gam_k)$) is equal to
$2|\lam-\mu|$. In particular, $\Bun_{G,\mu}^{\lam}$ is empty unless $\lam\geq\mu$.
\eth
\reft{dimension} is proved in \refss{dim}. Its proof is in fact unconditional (i.e. it does not assume the validity
of \refco{connect}), so by ``dimension'' we actually mean
``the dimension of every connected component
of $\Bun_{G,\mu}^{\lam}(\AA^2/\Gam_k)$''.

\medskip
\noindent
Another motivation for \refe{deflammu} is given by \reft{largekmain}.
\ssec{extreme}{Factorization and extremal connected components}
Here we are going to give another (perharps more intuitive) definition of $\Bun_{G,\mu}^{\lam}(\AA^2/\Gam_k)$.
The factorization morphism $\pi^a:\Bun_G^a(\AA^2)\to\BA^{(a)}$ is
$\Gamma_k$-equivariant (with respect to the action $z\mapsto\zeta z$ of
$\Gamma_k$ on $\BA^1$), and gives rise to the factorization morphism
$$
\pi_{\bmu}^{\blambda,a}:\Bun_{G,\bmu}^{\blambda,a}(\AA^2/\Gamma_k)\to(\AA^{(a)})^{\Gamma_k}\simeq
\BBA^{(b)},
$$
where $b=[a/k]$ is the integral part of $a/k$, and
$\BBA^1=\BA^1/\Gamma_k$ is also isomorphic to the affine line. However,
there is no a priori reason for $\pi_{\bmu}^{\blambda,a}$ to be surjective.
In fact, assuming \refco{connect} now we have the following:
\lem{minimalctwo}
Let $\blambda$ and $\bmu$ and $a$ be such that $\Bun_{G,\bmu}^{\blambda,a}$ is non-empty. Then
\begin{enumerate}
\item
There exists unique $d\in\NN$ such that $\frac{a-d}{k}\in\NN$ and such that the image of
$\pi_{\bmu}^{\blambda,d}$ consists of one point.
\item
In the above case let us set $c=\frac{a-d}{k}$. Then the image of $\pi_{\bmu}^{\blambda,a}$
is isomorphic to $\BBA^{(c)}$ which is embedded into $\BBA^{(b)}$ by adding a
multiple of $0\in\BBA^1$.
\end{enumerate}
\elem
\prf
Choose some $c$ such that $a-kc\geq 0$ and embed $\BBA^{(c)}$ into $\BBA^{(b)}$ in the way described above.
Assume now that there exists a divisor $D\in\BBA^{(c)}$
{\em disjoint from} $0\in\BBA^1$ which lies in the image of $\pi_{\bmu}^{\lambda,a}$.

Then
\eq{factor-bred}
(\pi_{\bmu}^{\blambda,a})^{-1}(D)\simeq
(\pi_{\bmu}^{\blambda,a-kc})^{-1}(0)\times(\pi^c)^{-1}(D)
\end{equation}
(here $\pi^c:\Bun_G^c(\AA^2)\to\BA^{(c)}$is the usual factorization
morphism, and we somewhat abusively identify $\BA^{(c)}$ and
$\BBA^{(c)}$). Hence if we choose $c$ to be the maximal number with
this property (i.e. that $D$ as above exists) and if we set
$d=a-kc$, then the image of the factorization morphism
$\pi_{\bmu}^{\blambda,d}$ is just one point 0. The uniqueness of $d$
as well as the second assertion of \refl{minimalctwo} follows
immediately from the factorization isomorphism \refe{factor-bred}.
\epr Loosely speaking, using the factorization isomorphism
\refe{factor-bred} the variety
$\Bun_{G,\bmu}^{\blambda,a}(\AA^2/\Gamma_k)$ can be generically
expressed in terms of $\Bun_{G,\bmu}^{\blambda,d}(\AA^2/\Gamma_k)$,
and of the usual $\Bun_G(\AA^2)$.

Thus for each couple $\bmu,\blambda\in\Lam^+_k$ there exist certain
{\em extremal components}
$\Bun_{G,\bmu}^{\blambda,d}(\AA^2/\Gamma_k)$ such that the image of
factorization morphism is just one point, and any other component
with the same $\bmu$ and $\blambda$ is of the kind
$\Bun_{G,\bmu}^{\blambda,d+kc}(\AA^2/\Gamma_k)$ for some $c\in\BN$.

\conj{connec} (a) For any $\bmu,\blambda\in\Lam_k^+$ there is a {\em
unique} extremal component
$\Bun_{G,\bmu}^{\blambda,d}(\AA^2/\Gamma_k)$ to be denoted
$^{\min}\Bun_{G,\bmu}^{\blambda}(\AA^2/\Gamma_k)$.

(b) There is a lift $\bmu_\infty:\ \GG_m\to G$
(resp. $\blambda_\infty:\ \GG_m\to G$) of
$\bmu:\ \Gamma_k\to G$ (resp. $\blambda:\ \Gamma_k\to G$) such that
$^{\min}\Bun_{G,\bmu}^{\blambda}(\AA^2/\Gamma_k)=
\Bun_{G,\bmu_\infty}^{\blambda_\infty}(\AA^2/\GG_m)$.
\footnote{The definition of $\Bun_G(\AA^2/\GG_m)$ is essentially the same as
that of $\Bun_G(\AA^2/\Gam_k)$; more details can be found in \refs{large}.}
\econj
\ssec{}{Remark} We will see in \refs{large} that
$\Bun_{\bmu_\infty}^{\blambda_\infty}(\AA^2/\GG_m)$ is
isomorphic to $\calW_{G,\bmu_{\infty}}^{\blambda_{\infty}}$.

\medskip
\noindent
In the rest of this Section we shall assume \refco{connec} as well. Note that as before
\refco{connec} is obvious for $k$=1 (by \cite{bfg}) and we shall see that it also holds
true for $G=\SL(n)$ in \refs{sln}.

\ssec{dwi}{The main conjecture}
Let $\blambda,\bmu\in\Lam^+_k$ and choose a
lift $\lambda=(k,\blambda,l)$ of $\blambda$ to a dominant weight of $G_{\aff}^{\vee}$.
Note that any other choice $\lambda'$ equals $\lambda+n\delta$ where $\delta=(0,0,1)$
is the minimal imaginary root of $G_{\aff}^{\vee}$, and $n\in\BZ$.
Among all the possible lifts $\mu$ of $\bmu$ there is a unique
{\em maximal} one $\mu_0$ such that $\mu_0$ has a nonzero
multiplicity in the integrable module $L(\lambda)$ with highest weight
$\lambda$, while $\mu_0+\delta$ has zero multiplicity. The set of
lifts $\mu_0-\BN\delta$ is called a {\em string} in the
terminology of integrable modules. To identify it with a corresponding
factorization string of connected components of
$\Bun_{G,\bmu}^{\blambda}(\AA^2/\Gamma_k)$, we will denote
$\Bun_{G,\bmu}^{\blambda,d+kc}(\AA^2/\Gamma_k)$
\footnote{Here, as before, we have $\Bun_{G,\bmu}^{\blambda,d}(\AA^2/\Gamma_k)=
{^{\min}\Bun}_{G,\bmu}^{\blambda}(\AA^2/\Gamma_k)$.}
by
$'\Bun^{\lambda}_{G,\mu_0-c\delta}(\AA^2/\Gamma_k)$ provisionally.
So for $\mu$ in the $\bmu$-string of $L(\lambda)$ we shall denote by
$'\Bun^{\lambda}_{G,\mu}(\AA^2/\Gamma_k)$ the corresponding connected
component. We shall also denote by $'\calU^{\lambda}_{G,\mu}(\AA^2/\Gamma_k)$ its
closure in the Uhlenbeck space of $\AA^2$. Note that
$'\calU^{\lambda+l\delta}_{G,\mu+l\delta}(\AA^2/\Gamma_k)=\
{}'\calU^{\lambda}_{G,\mu}(\AA^2/\Gamma_k)$ for arbitrary $l\in\BZ$.

The main claim of this paper is this:

\medskip
$\bullet\ $ One can (and should) think about the varieties
$\calU^{\lambda}_{G,\mu}(\AA^2/\Gamma_k)$ as analogs of the transversal slices
$\ocalW^{\lam}_{G,\mu}$ we considered in \refss{trans-finite}; similarly one should think about
$\Bun_{G,\mu}^{\lambda}$'s as affine analogs of the varieties  $\calW^{\lam}_{G,\mu}$.

\medskip
\noindent
Below is the main conjecture of this paper, which partially justifies the above principle. Some other (perhaps
more convincing) evidence may also be found in~\refs{large}.
\conj{Greenland}
\begin{enumerate}
\item
 For $\mu=(k,\bmu,m),\ \lambda=(k,\blambda,l)$ we have
$'\Bun^{\lambda}_{G,\mu}(\AA^2/\Gamma_k)=
\Bun^{\lambda}_{G,\mu}(\AA^2/\Gamma_k):=\Bun^{\blambda,a}_{G,\bmu}(\AA^2/\Gamma_k)$
where $a=k(l-m)+\frac{(\blambda,\blambda)}{2}-\frac{(\bmu,\bmu)}{2}$.
Hence $'\calU^{\lambda}_{G,\mu}(\AA^2/\Gamma_k)=
\calU^{\lambda}_{G,\mu}(\AA^2/\Gamma_k)$.
\item
We have
$\calU^\nu_{G,\mu}(\AA^2/\Gamma_k)\subset
\calU^\lambda_{G,\mu}(\AA^2/\Gamma_k)$ if and only if  $\lambda\geq\nu$ (that is, $\lambda-\nu$
is a positive linear combination of simple roots of $\grg_{\aff}^{\vee}$).
The corresponding closed embedding will be called ``adding defect
$(\lambda-\nu)\cdot(0,0)$'' (note that this part of the conjecture is compatible with \reft{dimension}).

\item
Let $V^{\lam}_{\mu}$ denote
the stalk of $\IC(\calU_{G,\mu}^{\lambda}(\AA^2/\Gamma_k))$ at the unique
torus-fixed point $pt$. This is a graded vector space. Then $V^{\lam}_{\mu}$ is concentrated in even degrees
and we have a canonical isomorphism
$$
(V^{\lam}_{\mu})_{-2i}\simeq \gr^F_iL(\lambda)_{\mu}
$$
(see~\refss{referee} for the definition of filtration $F$).
\end{enumerate}
\econj

Note that $\calU_{G,\mu}^{\lam}$ admits a $\CC^*$-action which contracts it to $pt$. Thus
the global intersection cohomology $\IH^*(\calU_{G,\mu}^{\lam})$ is also equal to $V^{\lam}_{\mu}$.

\ssec{aroundkostaff}{The $q$-analog of the Kostant multiplicity formula in the affine case}
We now want to describe the affine analog of \refss{aroundkost}.
In the same way as in \refss{aroundkost} let us introduce the generating functions:
$$
\IC^{\lam}_{\mu}(q):=\sum\limits_{i}\dim(V^{\lam}_{\mu})_{-2i}q^i;\quad
^eC^{\lam}_{\mu}(q):=\sum\limits_{i} \dim\gr_i^F L(\lam)_{\mu} q^i.
$$
Then the last assertion of \refco{Greenland} says that
$\IC^{\lam}_{\mu}(q)={^eC^{\lam}_{\mu}}$. We shall refer to this function as
{\em the $q$-character} of $L(\lam)_{\mu}$.

The $q$-character of $L(\lam)_{\mu}$ can be given another purely combinatorial definition.
Namely, let $\Lam_{\aff}^{\pos}$ denote the subsemigroup of $\Lam_{\aff}$ consisting of elements of
the form $\sum a_j \alp_j$ where $a_j\in \ZZ_{\geq 0}$ and $\alp_j$ are  positive coroots of $G_{\aff}$
(which is the same as positive roots of $G_{\aff}^{\vee}$). Note that $\Lam_{\aff}^{\pos}\subset \Lam_{\aff,0}$.
Let also $R_{\aff}$ denote the set of all roots of $G_{\aff}^{\vee}$ and let $R_{\aff}^+$
be the corresponding set of positive roots.
For all $\beta\in \Lam^{\pos}$ let us
introduce the Kostant partition function of $K_{\beta}(q)$ of $\beta$ by:
$$
\sum_{\beta\in\Lam^{\pos}_{\aff}}K_\beta(q):=
\prod_{\alpha\in R_{\aff}^+}(1-qe^\alpha)^{-1},
$$
Now we define following~\cite{Vi}
\eq{kostantaff}
C^\lambda_\mu(q):=
\sum_{w\in W_{\aff}}(-1)^{\ell(w)}K_{w\cdot\lambda-\mu}(q).
\end{equation}

Then we have the following
\conj{lusztig-brylinski-aff}For every $\lam,\mu\in\Lam_{\aff}^+$
We have $\IC^{\lam}_{\mu}(q)={^eC^{\lam}_{\mu}(q)}=C^{\lam}_{\mu}(q)$.
\econj

Note that the equality ${^eC^{\lam}_{\mu}(q)}=C^{\lam}_{\mu}(q)$ is purely representation-theoretic,
i.e. it doesn't involve any geometry in it.
This equality is the natural affine generalization of the main result of \cite{Bry}; in some special cases it was proved
in \cite{fgt} (in fact, I.~Grojnowski told us that the methods of \cite{fgt} allow to prove the equality
${^eC^{\lam}_{\mu}(q)}=C^{\lam}_{\mu}(q)$ in general but we were not able to reconstruct that proof).
It is also natural to expect that with the appropriate modifications
(as in Theorem 2.9 in \cite{joseph}) the equality
${^eC^{\lam}_{\mu}(q)}=C^{\lam}_{\mu}(q)$ can be generalized to the case when $\mu$ is not necessarily dominant.

\ssec{}{The full slices}In \refss{trans-finite} we introduced also the infinite-dimensional varieties $\calW_{\mu}$
(in that case $\mu$ was an element of $\Lam$, not of $\Lam_{\aff}$); more precisely, $\calW_{\mu}$ was an ind-scheme
and it might be thought of as a transversal slice to $\Gr_G^{\mu}$ in $\Gr_G$. It turns out that
we may construct and analog of $\calW_{G,\mu}$ in the affine case.

Namely, we consider the
inductive limit of the system of embeddings
$\calU^{\lambda}_{G,\mu}(\AA^2/\Gamma_k)
\hookrightarrow \calU^{\lambda+c\delta}_\mu(\AA^2/\Gamma_k)$ (for $c\in\BN$;
the embedding is just twisting by a multiple of $(0,0)\in\BBA^2$), and
denote it by $\calU_{G,\mu}(\AA^2/\Gamma_k)$. Note that the limit does not
depend on $\lambda$, as suggested by the notation (for various choices of $\lambda$
we obtain cofinal systems).

\sec{large}{The case of $k\ggg0$}
\ssec{}{}We want to understand what happens when we fix $\olam,\omu$ and $a$ and make $k$ very large.
Let us choose some $\bmu:\GG_m\to T\subset G$. As before it defines an action of $\GG_m$ on $\Bun_G^a(\AA^2)$.
Thus for each $\blambda,\bmu\in\Lam^+$ and $a\geq 0$ we may speak of $\Bun^{\blambda,a}_{G,\bmu}(\AA^2/\GG_m)$ and we
let $\calU_{G,\bmu}^{\blambda,a}(\AA^2/\GG_m)$ denote its closure in $\calU_G^a(\AA^2)$. It is obvious that the factorization
mapping $\pi_{\bmu}^{\blambda,a}$ sends all of $\calU_{G,\bmu}^{\blambda,a}(\AA^2/\GG_m)$ to the divisor $a\cdot 0$.

On the other hand,
it is clear that for $k\ggg0$ we have
$(\Bun_G^a(\AA^2))^{\Gam_k}=(\Bun_G^a)^{\GG_m}$. Hence if $k$ is large we have
$\Bun^{\blambda,a}_{G,\bmu}(\AA^2/\Gam_k)=\Bun^{\blambda,a}_{G,\bmu}(\AA^2/\GG_m)$.
\th{largekmain}
\begin{enumerate}
\item
The variety $\Bun_{G,\bmu}^{\blambda,a}(\AA^2/\GG_m)$
is empty unless $\lam\geq \mu$ and
$a=\frac{(\blambda,\blambda)}{2}-\frac{(\bmu,\bmu)}{2}$.
\item
In the above case (i.e. when $a=\frac{(\blambda,\blambda)}{2}-\frac{(\bmu,\bmu)}{2}$) there are canonical isomorphisms
$$
\Bun_{G,\bmu}^{\blambda,a}(\AA^2/\GG_m)\simeq \calW_{G,\bmu}^{\blambda};\qquad
\calU_{G,\bmu}^{\blambda,a}(\AA^2/\GG_m)\simeq
\overline{\calW}{}_{G,\bmu}^{\blambda}.
$$
\item
Assume that $\blambda\geq \bmu$ and $\blambda,\bmu\in\Lam^+_k$ for some $k>0$.
Let $a=\frac{(\blambda,\blambda)}{2}-\frac{(\bmu,\bmu)}{2}$. Assume that $k>a$.
Then $\Bun_{G,\bmu}^{\blambda,b}(\AA^2/\GG_m)$ is empty if $b<a$.
\end{enumerate}
\eth
\prf
Let us first concentrate on the first two assertions of \reft{largekmain}.
Let $\Gr_{G,BD}^n$ denote the Beilinson-Drinfeld Grassmannian classifying the following data:

\smallskip
a) An effective divisor $D=\sum a_i x_i$ on $\AA^1$ of degree $n$.

b) A $G$-bundle $\calF$ on $\PP^1$ trivialized away from the support of $D$.

\smallskip
\noindent
It is well-known that $\Gr_{G,BD}^n$ has a natural structure of an ind-scheme.
Note that we have a natural embedding $\Gr_G\hookrightarrow \Gr_{G,BD}^n$; the image of this embedding
consists of all the points of $\Gr_{G,BD}^n$ for which all $x_i$ are equal to zero.

Let $\bfC$ denote another copy of $\PP^1$; we denote by $\infty_{\bfC}$ and $0_{\bfC}$ the corresponding
0 and $\infty$ points.
Recall from \cite{bfg} that $\Bun_G^a(\AA^2)$ can be identified with the space $\Maps^a(\bfC,\Gr_{G,BD}^a)^0$
of all maps from
$\bfC\simeq\PP^1$ to $\Gr_{G,BD}^a$  which

\medskip
1) Have degree $a$ in the appropriate sense;

2) Send $\infty_{\bfC}$ to a point of $\Gr_{G,BD}^a$ corresponding to $\calF$ being trivial bundle with
the standard trivialization
(with $D$ being arbitrary).

\medskip
\noindent


Assume now that we are given a point of
$\Bun_{G,\bmu}^{\blambda,a}(\AA^2/\GG_m)$.
Then it is clear that the map $\pi^a$ sends this point to the
divisor $a\cdot 0$; also, the divisor $D$ must necessarily be of the
form $a\cdot 0_{\bfC}$,
that means that this points corresponds to a
mapping $f:\bfC\to \Gr_G$. The map $f$ must satisfy the
following properties:

$\bullet$ $f$ is of degree $a$

$\bullet$ $f(\infty_{\bfC})$ is the unique $G(\calO)$-invariant point of $\Gr_G$.

$\bullet$
$f(\tau z)=\bmu(\tau^{-1})(f(z))^{\tau^{-1}}$
where the super-script $\tau^{-1}$ in the right hand side stands for the action of $\tau^{-1}$ by
by loop rotation on $\Gr_G$.

$\bullet$ $s^{\bmu}f(0_{\bfC})$ lies in $(\Gr_G^{\blambda})^{\CC^*}$ (this follows from \refl{grastorsor}) .

It is clear that such an $f$ is uniquely
determined by its value at $1\in\bfC$. Let us consider the map $g=s^{\bmu}\cdot f$.
Then $g$ is again a map of degree $a$ which satisfies

a) $g(\tau z)=(g(z))^{\tau^{-1}}$,

b) $g(\infty_{\bfC})=s^{\bmu}$,

c) $g(0_{\bfC})\in (\Gr_G^{\blambda})^{\CC^*}$.

\noindent
The map $g$ again is uniquely determined by its value at 1.
In other words, evaluation of $g$ at 1 gives rise to an embedding $\eta$
of $\Bun_{G,\bmu}^{\blambda,a}(\AA^2/\GG_m)$ into $\Gr_G$.
\lem{breddd}
\begin{enumerate}
\item
The space of maps $g$ as above is empty unless $a=\frac{(\blambda,\blambda)}{2}-\frac{(\bmu,\bmu)}{2}$.
\item
In the latter case the embedding
$\eta$ is an isomorphism between $\Bun_{G,\bmu}^{\blambda,a}(\AA^2/\Gam_k)$ and $\calW_{G,\bmu}^{\blambda}$.
\end{enumerate}
\elem

\prf
It is well-known (cf. \cite{Fal}) that $\Pic(\Gr_G)\simeq \ZZ$;
different line bundles on $\Gr_G$
correspond canonically to different integral even symmetric invariant bilinear
forms on $\grg$. Let $\calL$ denote the generator
of $\Pic(\Gr_G)$ corresponding to our form $(\cdot,\cdot)$.
Then $\deg(g)$ is by definition the degree of the line bundle $g^*\calL$. We claim that for any map $g$ satisfying
the properties a),~b),~c) above one has
$\deg(g^*\calL)=\frac{(\blambda,\blambda)}{2}-\frac{(\bmu,\bmu)}{2}$. Indeed, it follows from
a) that the line bundle $g^*\calL$ on $\bfC\simeq \PP^1$ is $\GG_m$-equivariant. It also follows from b) and c)
that the action of $\GG_m$ in the fiber of $g^*\calL$ at $\infty_{\bfC}$ is given by
$\tau\mapsto\tau^{\frac{(\blambda,\blambda)}{2}}$
and the action of $\GG_m$ in the fiber of $g^*\calL$ at $0_{\bfC}$ is
given by $\tau\mapsto \tau^{\frac{(\bmu,\bmu)}{2}}$.
It is easy to see (according to the standard classification
of $\GG_m$-equivariant line bundles on $\PP^1$) that this implies that
$\deg (g^*\calL)=\frac{(\blambda,\blambda)}{2}-\frac{(\bmu,\bmu)}{2}$.
This proves the first assertion of \refl{breddd}.

To prove the second assertion, let us note that a point $x$ in $\Gr_G$ lies in $\calW^{\blambda}_{G,\bmu}$
if and only the following two properties are satisfied:
\eq{one}
\lim\limits_{\tau\to 0} x^{\tau}\in \Gr_G^{\blambda}.
\end{equation}
and
\eq{two}
\lim\limits_{\tau\to \infty} x^{\tau}=s^{\bmu}.
\end{equation}
Let us now prove that $\eta$ lands in $\calW_{G,\bmu}^{\blambda}$.
Indeed, let $g$ be as above and set $x=g(1)$. Then a),b) and c) above
imply that $x$ satisfies \refe{one} and \refe{two} which proves what we need.

On the other hand let us show that the map $\eta:\Bun_{G,\bmu}^{\blambda,a}\to \calW_{G,\bmu}^{\blambda}$
is surjective. Take any $x\in \calW_{G,\bmu}^{\blambda}$ and set
$g(z)=x^{z^{-1}}$ for any $z\in \CC^*$. Since $\Gr_G$ is ind-proper, it follows that
$g$ extends to a map $\PP^1\to\Gr_G$ which automatically satisfies a). Now the equations
\refe{one} and \refe{two} imply that $g$ also satisfies b) and c).
\epr

The first assertion of \refl{breddd} proves \reft{largekmain}(1). To prove the second assertion of
\reft{largekmain} we must show that $\eta$ extends to an isomorphism
between $\calU_{G,\bmu}^{\blambda,a}(\AA^2/\GG_m)$ and
$\ocalW^{\blambda}_{G,\bmu}$. Actually, the definition of $\eta$ extends word-by-word to
$\calU_{G,\bmu}^{\blambda,a}$; it is also easy to show in the same way as before
that in this way we get a map $\oeta:\calU_{G,\bmu}^{\blambda,a}\to
\ocalW_{G,\bmu}^{\blambda}$. We now must prove that this map is an isomorphism.

However, the
stratification \refe{stratification} of $\calU_G^a(\AA^2)$  implies that
$$
\calU_{G,\bmu}^{\blambda,a}(\AA^2/\GG_m)\simeq
\bigsqcup\limits_{\bmu\leq\blambda'\leq\blambda}\Bun_{G,\bmu}^{\blambda',
\frac{(\blambda',\blambda')}{2}-\frac{(\bmu,\bmu)}{2}}(\AA^2/\GG_m).
$$
Similarly,
$$
\ocalW^{\blambda}_{G,\bmu}=\bigsqcup_{\bmu\leq\blambda'\leq\blambda}
\calW_{G,\bmu}^{\blambda'}.
$$
This implies that $\oeta$ is a bijection on closed points.
Note that the union of the boundary
strata in $\ocalW^{\blambda}_{G,\bmu}$ is of codimension at least~2; hence $\oeta$
is an isomorphism away from codimension~2. On the other hand,
$\ocalW^{\blambda}_{G,\bmu}$ is known to be normal (cf. \cite{Fal}
and references therein). This implies that $\oeta^{-1}$ is defined
everywhere and since $\oeta$ is a bijection on closed points it follows that
$\oeta$ is is an isomorphism.

The proof of the third assertion of \reft{largekmain} is a word-by-word repetition
of the proof of the first assertion of \refl{breddd}.
\epr
\ssec{}{The case of large $k$}
Assume that we have $\blambda\geq \bmu$ such that
$\blambda,\bmu\in\Lam_{k}^+$. Then \reft{largekmain} produces certain
connected component of $\Bun_{G,\bmu}^{\blambda,a}(\AA^2/\Gam_k)$ for $a=\frac{(\blambda,\blambda)}{2}-\frac{(\bmu,\bmu)}{2}$ --- namely,
the one which contains $\Bun_{G,\bmu}^{\blambda,a}$. Moreover, it follows from \reft{largekmain} that for
$\frac{(\blambda,\blambda)}{2}-\frac{(\bmu,\bmu)}{2}<k$ this component is extremal.
 We don't know how to prove that
there are no other extremal components in this case; however, it follows from \reft{largekmain} that
in this case we have $^{\min}\Bun_{G,\bmu}^{\blambda}(\AA^2/\Gam_k)
=\Bun_{G,\bmu}^{\blambda,a}(\AA^2/\Gam_k)$.
This proves \refco{Greenland}(1) under the above assumptions. Indeed,
let us choose as before a lift $\lam=(k,\blambda,m)$ of $\blambda$ to $\Lam_{\aff,k}^+$.
Recall that in \refss{dwi} we introduced a weight $\mu_0\in\Lam_{\aff,k}^+$, which is a lift
of $\bmu$ to $\Lam_{\aff,k}^+$, and such that $L(\lam)_{\mu_0}\neq 0$ and $L(\lam)_{\mu_0+\del}=0$.
Then \refco{Greenland}(1) in this case is equivalent to saying that  $\mu=(k,\bmu,m)$.
To prove this let us note that the integrable $\grg_{\aff}^{\vee}$-module $L(\lam)$ contains
the finite-dimensional irreducible $\grg^{\vee}$-module $L(\blambda)$. It is equal to the direct sum of
all the $L(\lam)_{\nu}$ where $\nu$ is of the form $(k,\onu,m)$. Thus, since for a dominant weight
$\bmu$ satisfying $\bmu\leq \blambda$ we have $L(\blambda)_{\bmu}\neq 0$, it follows that
$L(\lam)_{(k,\bmu,m)}\neq 0$. On the other hand, we have $L(\lam)_{(k,\bmu,m)+\del}=0$
since $(k,\bmu,m)+\del=(k,\bmu,m+1)\not\leq (k,\blambda,m)=\lam$.

Even though we can't prove the absence of other extremal components in this case, we can now set
$\Bun_{G,\bmu}^{\blambda,a+kc}(\AA^2/\Gam_k)$ to be the {\em connected subvariety} of $\Bun_G(\AA^2/\Gam_k)$
obtained by adding defect $c\cdot (0,0)$ to the above extremal component and just ignore the other components
of $\Bun_G(\AA^2/\Gam_k)$ if they exist (of course, we believe that they don't exist). In this way, the notation
$\Bun_{G,\mu}^{\lam}$ acquires definite meaning, which is compatible with \refco{Greenland}(1).

Also, in this case we can now prove the dimension formula, claimed in the first assertion of \refco{Greenland}(2):
it follows immediately from \refl{gras-orbits}(3).
\ssec{}{Proof of the conjectures for large $k$}
Fix as before $\blambda$ and $\bmu$ such that $\blambda\geq\bmu$. According to \reft{largekmain}
the only component of $\Bun_{G,\bmu}^{\blambda}(\AA^2/\Gam_k)$ which survives
in the limit $k\to\infty$ is $\Bun_{G,\bmu}^{\blambda,a}(\AA^2/\Gam_k)$ where $a=\frac{(\blambda,\blambda)}{2}-\frac{(\bmu,\bmu)}{2}$.

It follows that in the limit $k\to \infty$
all the conjectures of \refs{formulation} hold in this case, except, possibly \refco{Greenland} whose
proof is explained below.

It is enough for us to assume that $\blambda\geq\bmu$ such that $\blambda,\bmu\in\Lam^+_k$ and
$\lam=(k,\blambda,m),\mu=(k,\bmu,m)$. According to the previous subsection, in this case
we have
\eq{fin-aff-mod}
L(\lam)_{\mu}=L(\blambda)_{\bmu}.
\end{equation}
Thus, using the assertion of
\reft{largekmain}, it follows that in order to prove \refco{Greenland} in this case, it is enough
to check that the filtration on the LHS of \refe{fin-aff-mod} defined by the principal nilpotent
$e_{\aff}$ of $\grg_{\aff}^{\vee}$ coincides with the filtration defined on the RHS of \refe{fin-aff-mod}
by the principal nilpotent $e$ of $\grg^{\vee}$. This is obvious, since the difference $e_{\aff}-e$
acts trivially on $L(\blambda)\subset L(\lam)$.

\sec{k=1}{The case of $k=1$}
In the case $k=1$ all the conjectures of \refs{formulation} follow immediately
from \cite{bfg} except for \refco{Greenland}(3). In this Section we explain
that
\refco{Greenland}(3) can be reduced to an explicit statement in representation
theory, which most probably follows from \cite{fgt},~section~11.
In particular, in this case
we prove that the
total stalk $V^{\lam}_{\mu}$ of $\IC(\calU_{G,\mu}^{\lambda}(\AA^2/\Gamma_k))$ at
the canonical point of $\calU_{G,\mu}^{\lambda}(\AA^2/\Gamma_k)$ is isomorphic to {\em some}
graded version of $L(\lam)_{\mu}$ (i.e. we prove the equality
$\IC^{\lam}_{\mu}(1)={^eC^{\lam}_{\mu}(1)}$). We also prove the equality
$\IC^{\lam}_{\mu}(q)=C^{\lam}_{\mu}(q)$ in the $ADE$ case.
\ssec{}{The structure of $\grg_{\aff}^{\vee}$}
Recall that the Langlands dual Lie algebra $\grg_{\aff}^\vee$ (in
general, twisted) admits the following realization. There exist

1)  A simple Lie algebra $\grg'$,

2) Its outer automorphism $\sigma$ of order
$m=1,2,3$,

3) A primitive $m$-th root of unity $\omega$

\noindent
such that
$$
\grg_{\aff}^\vee\simeq\CC c\oplus \CC d \oplus
\bigoplus_{n\in{\mathbb Z}}\grg'_n\otimes t^n
$$
where $c$ is the
central element, $d$ is the loop rotations element, and
$\grg'_n:=\{x\in\grg':\ \sigma(x)=\omega^nx\}$. Let $\grt'
\subset\grg'$ be a $\sigma$-invariant Cartan subalgebra, and let
$\grt'_n:=\grt'\cap\grg'_n$. We set
$$
\gra^-:=\bigoplus_{n<0}\grt'_n\otimes t^n.
$$
This is naturally a ${\mathbb
  Z}^-$-graded vector space. The symmetric algebra
$\operatorname{Sym}\gra^-$ acquires the induced ${\mathbb
  Z}^{\leq0}$-grading,  $\operatorname{Sym}\gra^-=
\bigoplus_{n\in{\mathbb Z}^{\leq0}}(\operatorname{Sym}\gra^-)_n$.
\ssec{}{The basic integrable module over $\grg_{\aff}^{\vee}$}
Let $\lambda=(1,0,0)$ be the basic fundamental weight of
$\grg_{\aff}^\vee$. Note that this is the only integrable weight of
$G^\vee_{\aff}$ at level 1 since $G$ is simply connected.
Then the only dominant weights $\mu$ of $\grg_{\aff}^{\vee}$
for which $L(\lambda)_{\mu}\neq 0$
are of the form
$\mu=\lambda+n\delta=(1,0,n),\ n\in{\mathbb Z}^{\leq0}$.
According to~\cite{Lep} (it is essential here that $\grg_{\aff}^\vee$
is Langlands dual of an {\em untwisted} affine Lie algebra),
we have $L(\lambda)_{\lambda+n\delta}\simeq
(\operatorname{Sym}\gra^-)_n$ for any $n\in{\mathbb Z}^{\leq0}$.

Now note that $\grg'_0$ is another simple Lie algebra. We choose its
principal nilpotent element $e'$. According to Kostant's theorem,
$\dim(\grg'_m)^{e'}=\dim(\grt'_m)$. Hence we have a
(noncanonical) isomoprhism of graded vector spaces $\gra^-
\simeq\bigoplus_{m<0}(\grg'_m)^{e'}$. The symmetric algebra
$\operatorname{Sym}(\bigoplus_{m<0}(\grg'_m)^{e'})$
acquires the induced ${\mathbb Z}^{\leq0}$-grading,
$\operatorname{Sym}(\bigoplus_{m<0}(\grg'_m)^{e'})=:
\bigoplus_{n\in{\mathbb Z}^{\leq0}}\operatorname{Sym}_n$.

\ssec{}{Comparison with IC-stalks}
According to~\cite{bfg}, the total IC stalk $V^{\lam}_{\lam+n\del}$ (with grading disregarded)
of $\operatorname{IC}({\mathcal U}^\lambda_{G,\lambda+n\delta}({\mathbb
  A}^2))$ at the canonical point is isomorphic to $\operatorname{Sym}_n$. Thus we have proved
a noncanonical isomorphism between $V^{\lam}_{\lam+n\del}$ and $L(\lambda)_{\lambda+n\delta}$, for any
$n\in{\mathbb Z}^{\leq0}$.

On the other hand, the vector space $(\grg'_m)^{e'}$ also has another grading, coming from a Jacobson-Morozov
triple $(e',h',f')$ in $\grg'_0$
(cf. the end of Section 7 of \cite{bfg} for more details); this grading
extends to each $\Sym_n$ (making
$\operatorname{Sym}(\bigoplus_{m<0}(\grg'_m)^{e'})$
a bigraded vector space). According to Theorem 7.10 of \cite{bfg}
the second grading on $\Sym_n$ (suitably normalized) corresponds to the
cohomological grading on $V^{\lam}_{\lam+n\del}$. Comparing with~Corollary~2
of~\cite{Vi} we see that in case $\mathfrak g$ is of type $ADE$
the generating function of graded dimensions of
$V^\lambda_{\lambda+n\delta}$ coincides with the level one $t$-string function
$\alpha^\lambda_\lambda(t)$, and hence
${\operatorname{IC}}^\lambda_{\lambda+n\delta}(q)=C^\lambda_{\lambda+n\delta}(q)$.
Thus, in order to prove the full statement of \refco{lusztig-brylinski-aff}
in this case we must check that
$\Sym_n$ is isomorphic to $L(\lam)_{\lam+n\del}$ as {\em graded} vector space.
We don't know how to do this, but
according to C.~Teleman and I.~Grojnowski (private communication)
one can probably use the techniques of \cite{fgt} in order to construct such an isomorphism.
\sec{sln}{The case of $G=\SL(N)$}
In this Section we prove most of the conjectures of \refs{formulation} in the case $G=\SL(N)$. This Section should
be thought of as an affine analog of the paper \cite{MVy},
where the authors discuss the relations between the
(nonaffine) $A_{k-1}$ quiver varieties, the transversal slices in the
(classical) affine Grassmannians of $\SL(N)$, and Howe's skew
$(\GL(k),\GL(N))$-duality.
\ssec{}{Cyclic quiver varieties and their IC-sheaves}
We consider the cyclic quiver with $k$ vertices corresponding to the
characters of $\Gam_k$~\cite{nakajima}.
For two $\Gamma_k$-modules $V,W$ we denote by
$\mathfrak{M}_0(V,W)$ the affine (singular, in general) Nakajima quiver variety (cf. \cite{nakajima});
\footnote{Of course here the word ``affine" is used in the sense of algebraic geometry and not
in the sense of ``affine Lie algebras".}
we also
denote by $\mathfrak{M}_0^{\operatorname{reg}}(V,W)$ its regular part.
Let $N=\dim W\geq2,\ a=\dim V$. We set $G=\SL(N)$. Note that in this case $\grg_{\aff}\simeq\grg_{\aff}^{\vee}$.

Let us write $w=\underline{\dim}W=(w_1,\ldots,w_k)$ if

\medskip
1) The multiplicity of the trivial representation of $\Gam_k$ in $W$ is $w_k$.

2) For $i=1,...,k-1$ the multiplicity of the character $\zeta\mapsto\zeta^i$ in $W$ is $w_i$.

\medskip
Similarly, one may introduce $v=\underline{\dim}V=(v_1,\ldots,v_k)$.
Let $\omega_i$ be the $i$-th fundamental
weight of $\fgl(k)_{\aff}$; by definition it is of level 1. In the standard notation we write
$\omega_i=(1,E_1+\ldots+E_i,0)$. Let $\alpha_i$ be the $i$-th simple
root of $\fgl(k)_{\aff}$; in the standard notation
we have $\alpha_i=(0,E_i-E_{i+1},0)$ for $i=1,\ldots,k-1$, and
$\alpha_0=\alpha_k=(0,-E_1+E_k,1)$. Then Nakajima's criterion says that
$\mathfrak{M}_0^{\operatorname{reg}}(V,W)$ is nonempty if and only if the
following two properties hold:

$\bullet\ $ The $\fgl(k)_{\aff}$-weight
$\sum\limits_{i=1}^kw_i\omega_i-\sum\limits_{i=1}^kv_i\alpha_i$ is dominant.

$\bullet\ $ The multiplicity of $\sum\limits_{i=1}^kw_i\omega_i-\sum\limits_{i=1}^kv_i\alpha_i$
in the integrable $\fgl(k)_{\aff}$-module with highest weight
$\sum\limits_{i=1}^kw_i\omega_i$ is non-zero.

In this case Theorem~5.2 of \cite{nakajima} implies the following result:
\th{nakajimaic}
The stalk of $\IC(\mathfrak{M}_0(V,W))$ at the most singular point has the following properties:

a) It is concentrated in the even degrees.

b) Its total dimension (with grading disregarded) equals the multiplicity of the
integrable $\fgl(k)_{\aff}$-module
$L(\sum\limits_{i=1}^kw_i\omega_i-\sum\limits_{i=1}^kv_i\alpha_i)$ in the tensor
product $\bigotimes_{i=1}^kL(\omega_i)^{\otimes w_i}$ of
fundamental representations $L(\omega_i)$ of $\fgl(k)_{\aff}$.
\eth

\ssec{}{Interpretation via $\SL(N)$}
Let $G=\SL(N)$. Recall the bijection $\Psi_{N,k}$ of~\refss{vee} between
$\Lambda^+_k$
and the conjugacy classes of homomorphisms $\Gamma_k\to\SL(N)$, i.e.
isomorphism classes of $N$-dimensional $\Gamma_k$-modules of determinant 1.
We will write down the elements $\bmu\in\Lambda^+_k$ as generalized Young
diagrams, i.e. the sequences of integers $(\mu_1\geq\mu_2\geq\ldots\geq\mu_N)$
such that $\mu_1-\mu_N\leq k$, and $\mu_1+\ldots+\mu_N=0$.

\lem{bij}
If $\Psi_{N,k}(\bmu)=(w_1,\ldots,w_k)$, then $w_i$ equals the number of
entries among $\{\mu_1,\ldots,\mu_N\}$ congruent to $i$ modulo $k$,
for any $i\in\BZ/k\BZ$. \qed
\elem
If we view $(w_1,\ldots,w_k)$ as a dominant weight
$w_1\omega_1+\ldots+w_{k-1}\omega_{k-1}$ of $\operatorname{PSL}(k)$ then
there is a unique way to write it as a generalized Young diagram
$\bnu=(\nu_1\geq\nu_2\geq\ldots\geq\nu_k)$ such that $\nu_i-\nu_{i+1}=w_i$ for
any $1\leq i\leq k-1$, and $\nu_1-\nu_k\leq N$, and
$\nu_1+\ldots+\nu_k=0$. If $(w_1,\ldots,w_k)=\Psi_{N,k}(\bmu)$, we will write
$\bnu={}^t\bmu$ (transposed generalized Young diagram). For the corresponding
level $N$ weight of $\widehat{\mathfrak{gl}}(k)$ we have
$\sum_{i=1}^kw_i\omega_i=\sum_{j=1}^N\omega_{\mu_i}$. Here $\omega_{\mu_i}$
is understood as $\omega_{\mu_i\pmod{k}}$.


Similarly, let us consider the level $N$ dominant $\fgl(k)_{\aff}$-weight
$\sum\limits_{i=1}^kw_i\omega_i-\sum\limits_{i=1}^kv_i\alpha_i$.
Its projection to $\hatLam$ can be written as
$\sum\limits_{i=1}^kw'_i\omega_i$;
note that $\sum\limits_{i=1}^k w_i'=\sum\limits_{i=1}^k w_i=N$.
We now consider the corresponding generalized Young diagram $\blambda$ such
that $\Psi_{N,k}(\blambda)=(w_1',...,w_k')$ and we view $\blambda$ as a
dominant weight of $\widehat{\fsl(N)}$ at level $k$.
Let us now introduce $\mu,\lam\in\Lam_{\aff,k}^+$ by setting
$$
\mu=(k,~\bmu,~-\frac{a+\frac{(\bmu,\bmu)}{2}-\frac{(\blambda,\blambda)}{2}}{k}),
\qquad \lambda=(k,\blambda,0).
$$

Then we have the following
\lem{nak-uhl}
There exist natural isomorphisms
$\mathfrak{M}_0^{\operatorname{reg}}(V,W)\simeq \Bun_{G,\mu}^{\lam}$ and
$\mathfrak{M}_0(V,W)\simeq\calU_{G,\mu}^\lambda(\BA^2/\Gamma_k)$
where $\lam$ and $\mu$ are constructed as above.
\elem
\prf
The isomorphism $\mathfrak{M}_0^{\operatorname{reg}}(V,W)\simeq
\Bun_{G,\mu}^{\lam}$
is explained in \cite{nakajima}. The second isomorphism follows from the first and from
Theorem 5.12 of \cite{bfg}.
\epr

We now claim that \reft{nakajimaic} and \refl{nak-uhl} imply also that
$\dim V^{\lam}_{\mu}=
\dim L(\lam)_{\mu}$. In other words we can't prove the full statement
of~\refco{Greenland}(3) which
is essentially equivalent to the equality
$\IC^{\lam}_{\mu}(q)={^eC^{\lam}_{\mu}(q)}$, but we can show
that $\IC^{\lam}_{\mu}(1)={^eC^{\lam}_{\mu}(1)}$. From this,~\refco{Greenland}(1)
follows for free since the nonemptiness of $\Bun^\lambda_{G,\mu}$ is equivalent
to the nonvanishing of $\IC^\lambda_\mu(1)$, and hence the geometric and
representation-theoretic strings coincide. Moreover, the argument below shows
as well that the total dimension of
$\IC(\calU^\lambda_{\SL(N),\mu}(\BA^2/\Gamma_k))$-stalk at the generic point of
$\calU^\nu_{\SL(N),\mu}(\BA^2/\Gamma_k)$ equals $\dim L(\lambda)_\nu$, which
implies~\refco{Greenland}(2).

In order to verify $\dim V^{\lam}_{\mu}=
\dim L(\lam)_{\mu}$ it is enough to establish the following:
\prop{duality}
\eq{dual}
\Hom_{\fgl(k)_{\aff}}(L(\sum\limits_{i=1}^kw_i\omega_i-\sum\limits_{i=1}^kv_i\alpha_i),\bigotimes_{i=1}^kL(\omega_i)^{\otimes w_i})\simeq
L(\lambda)_\mu
\end{equation}
\eprop

\medskip
\noindent
{\bf Warning.} The notation in the formulation of \refp{duality} is a little bit misleading: the reader should
remember that the integrable modules appearing in the left hand side of \refe{dual} are modules over
$\fgl(k)_{\aff}$, and the module $L(\lam)$,
which appears in the right hand side, is a module over $\fsl(N)_{\aff}$
of level $k$.

\medskip
\prf
In effect, this is an instance of I.~Frenkel's level-rank duality
(see~\cite{Fr}). We will follow the notations of Appendix~A.1 of~\cite{Na}
(for a detailed exposition see~\cite{Ha}).
Inside $\fgl(kN)_{\aff}$ we consider the subalgebra
$\fgl(k)(N)$ (see \cite{Fr}, \S1.2, page 83) which is
$\hgl(k)+\hsl(N)$ with
the common central extension (and no loop rotation). The restriction of the
basic $\fgl(kN)_{\aff}$-module (that is, degree 0 part of the semiinfinite
wedge module) to this subalgebra has a simple spectrum
(\cite{Fr}, Proposition 2.2). Moreover, it is shown in~\cite{Fr} that:

\medskip
a)  As a $\hgl(k)$-module, this
basic module is a sum of tensor products of $N$-tuples of fundamental
$\hgl(k)$-modules ({\em loc. cit.}). This decomposition into direct sum
of tensor products corresponds to the decomposition into weight spaces
of the Cartan sublagebra of $\fsl(N)\subset \hsl(N) $ in such a way that
the module $\bigotimes_{i=1}^kL(\omega_i)^{\otimes w_i}$ over $\hgl(k)$
corresponds
to the weight $\bmu$ of $\fsl(N)$ where $\bmu$ and $w$ are related as above.

b) As a $\hgl(k)+\hsl(N)$ module the basic representation of
$\fgl(kN)_{\aff}$ is isomorphic to
the direct sum of modules of the form
$L(\sum\limits_{i=1}^kw_i\omega_i-\sum\limits_{i=1}^kv_i\alpha_i)\ten
L(\lam)$ where $v$ and $\lam$  are related as above
(here the first factor is a $\hgl(k)$-module and the second
factor is a $\hsl(N)$-module).

\medskip
This immediately proves \refp{duality} with the
loop rotation action disregarded. In other words, we get an isomorphism
\eq{dual'}
\Hom_{\hgl(k)}(L(\sum\limits_{i=1}^kw_i\omega_i-\sum\limits_{i=1}^k
v_i\alpha_i),\bigotimes_{i=1}^kL(\omega_i)^{\otimes w_i})\simeq
L(\lambda)_{\bmu}
\end{equation}
We need to check that the decomposition of both sides of \refe{dual'}
with respect to the loop rotation
gives rise to the isomorphism \refe{dual}.
This is nothing but the formulae~(A.6,A.7) of~\cite{Na}.
In effect, in notations of {\em loc. cit.}, we have to check that
\eq{nakaj}
\langle d^X,\lambda-\mu\rangle=
\frac{a+\frac{(\bmu,\bmu)}{2}-\frac{(\blambda,\blambda)}{2}}{k}
\end{equation}
According to~(A.7) of {\em loc. cit.}, we have
$\langle d^X,\lambda-\mu\rangle=v_k+\langle d,M(\blambda)\rangle-
\langle d,M(\bmu)\rangle$. Here $\langle d,M(\bmu)\rangle$ is the energy
level in the semiinfinite wedge $\fgl(kN)_{\aff}$-module of the highest vector
of the tensor
product $\bigotimes\limits_{p=1}^NL(\omega_{\mu_i})$ of the fundamental
$\fgl(k)_{\aff}$-modules.

\lem{negat}
$\langle d,M(\bmu)\rangle$ (resp. $\langle d,M(\blambda)\rangle$)
is the sum of negative entries of the corresponding
generalized Young diagram $\bmu$ (resp. $\blambda$). \qed
\elem

Adding a multiple of the imaginary root $\delta=(1,\ldots,1)$ to
$v=(v_1,\ldots,v_k)$, we reduce the proof of~\refe{nakaj} to the case
$v_k=0$. In this case $a=\langle\check\rho_{\SL(k)},
{}^t\bmu-{}^t\blambda\rangle$. Here $\check\rho_{\SL(k)}=\frac{1}{2}(k-1,k-3,
\ldots,3-k,1-k)$. Thus it remains to check $\langle\check\rho_{\SL(k)},
{}^t\blambda\rangle=-\frac{(\blambda,\blambda)}{2}-k\langle d,
M(\blambda)\rangle$. Recall that $^t\blambda$ is the sum of fundamental weights
$\omega_b$ of $\SL(k)$ (and $\omega_b$ is understood as $\omega_{b\pmod{k}}$).
If $^t\blambda_b=\omega_b$, and $b>0$, then
$\langle\check\rho_{\SL(k)},{}^t\blambda_b\rangle=\frac{1}{2}(bk-b^2)$,
and if $b<0$, then
$\langle\check\rho_{\SL(k)},{}^t\blambda_b\rangle=\frac{1}{2}(-bk-b^2)$.
Summing up these equalities over all $b$, we obtain
$\langle\check\rho_{\SL(k)},{}^t\blambda\rangle=-\frac{1}{2}\sum_bb^2
+\frac{1}{2}\sum_{b>0}bk-\frac{1}{2}\sum_{b<0}bk=
-\frac{(\blambda,\blambda)}{2}-k\sum_{b<0}b=
-\frac{(\blambda,\blambda)}{2}-k\langle d,M(\blambda)\rangle$ (the last
equality holds by~\refl{negat}). This completes the proof of~\refp{duality}.
\epr

\ssec{dim}{Dimensions}
In this section we prove \reft{dimension}.
We assume that ${\mathbf S}={\mathbb P}^2$, and pick up a $G$-bundle
${\mathcal F}\in\Bun^\lambda_{G,\mu}=\Bun^{\blambda,a}_{G,\bmu}$ in notations
of~\refd{uhl}. We have $\dim\Bun^\lambda_{G,\mu}=
\dim T_{\mathcal F}\Bun^\lambda_{G,\mu}$. Furthermore, we have
$T_{\mathcal F}\Bun^\lambda_{G,\mu}=
H^1({\mathbb P}^2,\opn{ad}_{\mathcal F}(-1))^{\Gamma_k}$
where $\opn{ad}_{\mathcal F}$ stands for the $\Gamma_k$-equivariant vector
bundle associated to $\mathcal F$ and the adjoint representation of $G$.
Now $\opn{ad}_{\mathcal F}$ can be described in the ADHM terms of the cyclic
$A_{k-1}$-quiver, and it has certain dimension vectors $v=(v_1,\ldots,v_k)$
and $w=(w_1,\ldots,w_k)$ attached to it. We have
$\dim H^1({\mathbb P}^2,\opn{ad}_{\mathcal F}(-1))^{\Gamma_k}=v_k$. So we must
prove $2|\lambda-\mu|=v_k$.

Recall that $2|\lambda-\mu|=\langle2\check\rho_G,\blambda-\bmu\rangle+
2\check{h}(l-m)=\langle2\check\rho_G,\blambda-\bmu\rangle+\frac{2\check{h}}{k}
(a+\frac{(\bmu,\bmu)}{2}-\frac{(\blambda,\blambda)}{2})$ where $\check{h}$ is
the dual Coxeter number of $G$. Note that $\opn{ad}_{\mathcal F}\in
\Bun^{\opn{ad}_*\blambda,2\check{h}a}_{\SL(\mathfrak{g}),\opn{ad}_*\bmu}$ where
$\opn{ad}_*\blambda$ and $\opn{ad}_*\bmu$ are the coweights (in fact, coroots)
of $\SL(\mathfrak{g})$. We have $(\opn{ad}_*\blambda,\opn{ad}_*\blambda)=
2\check{h}(\blambda,\blambda)$, and $(\opn{ad}_*\bmu,\opn{ad}_*\bmu)=
2\check{h}(\bmu,\bmu)$ since the Dynkin index of the adjoint representation
of $G$ is $2\check{h}$. Furthermore, if we write down $\opn{ad}_*\blambda$
(resp. $\opn{ad}_*\bmu$) as a generalized Young diagram (with $\dim\mathfrak g$
entries summing up to zero) then $\langle2\check\rho_G,\blambda\rangle$
(resp. $\langle2\check\rho_G,\bmu\rangle$) equals the sum of positive
entries of the corresponding Young diagram, which will be denoted
$-\langle d,M(\opn{ad}_*\blambda)\rangle$ (resp.
$-\langle d,M(\opn{ad}_*\bmu)\rangle$) in accordance with~\refl{negat}.

All in all, we must prove $$v_k=\langle d,M(\opn{ad}_*\bmu)\rangle
-\langle d,M(\opn{ad}_*\blambda)\rangle+
\frac{1}{k}\left(c_2(\opn{ad}_{\mathcal F})+
\frac{(\opn{ad}_*\bmu,\opn{ad}_*\bmu)}{2}-
\frac{(\opn{ad}_*\blambda,\opn{ad}_*\blambda)}{2}\right).$$
However, this is exactly
what was established in the proof of~\refp{duality}, compare~\refe{nakaj} and
the next line after it.
This completes the proof of \reft{dimension}.


\bigskip
\footnotesize{
{\bf A.B.}: Department of Mathematics, Brown University, 151 Thayer St., Providence RI
02912, USA;\\
{\tt braval@math.brown.edu}}

\footnotesize{
{\bf M.F.}: IMU, IITP and State University Higher School of Economics,\\
Department of Mathematics, 20 Myasnitskaya st, Moscow 101000 Russia;\\
{\tt fnklberg@gmail.com}}

\end{document}